\documentclass{amsart}
\usepackage{amssymb} 
\usepackage{color}
\usepackage[all]{xy}

 \newtheorem{theorem}{Theorem}[section]

 \newtheorem{proposition}[theorem]{Proposition}
 \newtheorem{lemma}[theorem]{Lemma}
 \newtheorem{corollary}[theorem]{Corollary}
 \newtheorem{remark}[theorem]{Remark}
 \newtheorem{definition}[theorem]{Definition}
 \newtheorem{notation}[theorem]{Notation}

 \numberwithin{equation}{section}

\def\cat{{\rm Cat}_{\infty}}
\def\wcat{\widehat{\rm Cat}_{\infty}}
\def\prl{{\rm Pr}^L}

\def\subdelta#1{({\mathbf\Delta}_{/[#1]})^{\rm op}}

\def\bialgebra{\mathrm{Bialg}}
\def\bicomod{\mathrm{Bicomod}}
\def\rcomod{\mathrm{Rcomod}}
\def\mons{\mathrm{Mon}_{\mathcal{O}}}
\def\monlax{\mathrm{Mon}_{\mathcal{O}}^{\rm lax}}
\def\monoplax{\mathrm{Mon}_{\mathcal{O}}^{\rm oplax}}

\begin{document}

\title
%[Multiplicative structures on comodules
%in higher categories]
{Multiplicative structures on comodules
in higher categories}
\author{Takeshi Torii}
\address{Department of Mathematics, 
%Faculty of Sience, 
Okayama University,
Okayama 700--8530, Japan}
\email{torii@math.okayama-u.ac.jp}
%\thanks{}

\subjclass[2020]{Primary 18N60; Secondary 18N70, 18M50, 55U40}
\keywords{Comodule, bialgebra,
  duoidal $\infty$-category, $\infty$-operad}

\date{March 3, 2025 (version~1.0)}

\begin{abstract}
In this paper
we study multiplicative structures
on comodules over bialgebras
in the setting of $\infty$-categories.
We show that the $\infty$-category
of comodules over an $(\mathcal{O},\mathbf{Ass})$-bialgebra
in a mixed $(\mathcal{O},\mathbf{Ass})$-duoidal
$\infty$-category has the structure
of an $\mathcal{O}$-monoidal $\infty$-category
for any $\infty$-operad $\mathcal{O}$.
\end{abstract}

\maketitle

\section{Introduction}

Comodules over coalgebras
are fundamental algebraic structures,
which are dualizations of modules over algebras.
If comodules are defined over a bialgebra,
then they obtain a multiplicative structure
by using the product on the bialgebra.
In this paper we investigate multiplicative structures
on comodules over bialgebras
in the setting of $\infty$-categories.

Comodules have many applications in various areas of mathematics
ranging over category theory, representation theory,
differential geometry, algebraic topology, algebraic geometry, 
non-commutative geometry and mathematical physics.
In algebraic topology
the Adams spectral sequence
is a useful tool to compute the stable homotopy groups
(\cite{Adams1, Adams2}).
Its $E_2$-page can be described
in terms of the comodule structure on mod $p$ homology
over the dual Steenrod algebra.
%Milnor (\cite{Milnor})
%has given a nice description
%of the structure of the dual Steenrod algebra.
The Adams-Novikov
spectral sequence
is a variant of Adams spectral sequence
based on the complex cobordism $MU$,
which predicts
nilpotence and periodicity 
in the stable homotopy category
(\cite{MRW, Morava, Ravenel}).
Its $E_2$-page can be described
in terms of the comodule structure on $MU$-homology
over the Hopf algebroid $MU_*(MU)$.

Motivated by these examples in stable homotopy theory
we have studied $\infty$-categories of
comodule spectra in \cite{Torii}.
Although we did not discuss it
%multiplicative structures on comodule spectra
in \cite{Torii},
we realized that it is necessary to 
develop the theory of multiplicative structures
on comodules in higher category theory.

In the classical setting
we can formulate categories of comodules
as follows.
Let $(\mathcal{M},\otimes)$
be the symmetric monoidal category
of $k$-modules for a commutative ring $k$.
For an algebra $A$ over $k$,
the category ${}_A\mathrm{BMod}_A(\mathcal{M})$
of $A$-$A$-bimodules over $k$
has the structure of a monoidal category
via the relative tensor product $\otimes_A$.
Furthermore,
the category $\mathrm{RMod}_A(\mathcal{M})$
of right $A$-modules is right tensored
over the monoidal category ${}_A\mathrm{BMod}_A(\mathcal{M})$.
By using this tensored structure,
we can set up a category $\rcomod_{(A,\Gamma)}$
of right $\Gamma$-comodules in $\mathrm{RMod}_A(\mathcal{M})$
for a coalgebra
$\Gamma$ in ${}_A\mathrm{BMod}_A(\mathcal{M})$.

Furthermore,
when $A$ is commutative,
the category ${}_A\mathrm{BMod}_A(\mathcal{M})$
has another tensor product,
which is given by
$M\boxtimes N=A\otimes_{A\otimes A} (M\otimes N)
\otimes_{A\otimes A}A$.
%for $A$-$A$-bimodules $M,N$.
Equipped with this tensor product $\boxtimes$ and
the relative tensor product $\otimes_A$,
the category
${}_A\mathrm{BMod}_A(\mathcal{M})$ has the structure
of a duoidal category (\cite{Aguiar-Mahajan}).

The notion of duoidal category was introduced
by Aguiar-Mahajan in \cite{Aguiar-Mahajan}
by name of $2$-monoidal category,
and the term of duoidal category
was proposed by Street in \cite{Street}.
A duoidal category has two tensor products
in which one is (op)lax monoidal with respect to
the other.
We can regard a duoidal category as a pseudomonoid
in the $2$-category of monoidal categories
and (op)lax monoidal functors.

When $A$ is commutative,
$\mathrm{RMod}_A(\mathcal{M})$
is also a monoidal category
by a similar tensor product to $\boxtimes$.
Furthermore,
$\mathrm{RMod}_A(\mathcal{M})$
is a right module over
the pseudomonoid ${}_A\mathrm{BMod}_A(\mathcal{M})$
in the $2$-category of monoidal categories
and lax monoidal functors.
We can consider a bialgebra
$\Gamma$ in the duoidal category
${}_A\mathrm{BMod}_A(\mathcal{M})$.
In this case
the category $\rcomod_{(A,\Gamma)}$
has the structure of a monoidal
category such that the forgetful
functor $\rcomod_{(A,\Gamma)}\to
\mathrm{RMod}_A(\mathcal{M})$
is strong monoidal.

In this paper we would like to
generalize the above multiplicative structure
on comodules in the setting of $\infty$-categories.
As in the classical setting,
for a good monoidal $\infty$-category $\mathcal{M}$,
we can consider an $\infty$-category
$\rcomod_{(A,\Gamma)}(\mathcal{M})$
of right $\Gamma$-comodules
in $\mathrm{RMod}_A(\mathcal{M})$
for an algebra object $A$ in $\mathcal{M}$
and a coalgebra object $\Gamma$
in ${}_A\mathrm{BMod}_A(\mathcal{M})$.

In order to study multiplicative structures
on comodules in the setting of $\infty$-categories,
we need a notion of duoidal $\infty$-category.
Duoidal $\infty$-categories have been introduced
in \cite{Torii1}
which are generalizations of duoidal categories
in the setting of $\infty$-categories.
Furthermore,
we generalized duoidal $\infty$-categories
to higher monoidal $\infty$-categories
in \cite{Torii2}.
As examples of duoidal $\infty$-categories,
we have shown that 
$\infty$-categories of operadic modules
have the structure of duoidal $\infty$-categories
under some circumstances
in \cite{Torii4}.
We have also given other examples
of duoidal $\infty$-categories
which are obtained as endomorphism
$\infty$-categories of map monoidales
in monoidal $(\infty,2)$-categories
in \cite{Torii5}.

We will investigate an $\mathcal{O}$-monoidal
structure on $\rcomod_{(A,\Gamma)}(\mathcal{M})$
for an $\infty$-operad $\mathcal{O}^{\otimes}$.
Suppose that $\mathcal{M}$ has in addition
an $\mathcal{O}$-monoidal
structure which is compatible with the monoidal structure.
In other words,
we suppose that
$\mathcal{M}$ is an $(\mathcal{O}\times\mathbf{Ass})$-monoidal
$\infty$-category.
When $A$ is an $(\mathcal{O}\times\mathbf{Ass})$-algebra
object of $\mathcal{M}$,
the $\infty$-category
${}_A\mathrm{BMod}_A(\mathcal{M})$
is a monoid object in the $\infty$-category of $\mathcal{O}$-monoidal
$\infty$-categories and lax $\mathcal{O}$-monoidal functors.
In other words,
${}_A\mathrm{BMod}_A(\mathcal{M})$
has the structure of
a mixed $(\mathcal{O},\mathbf{Ass})$-duoidal $\infty$-category.
In this case
we can consider an $(\mathcal{O},\mathbf{Ass})$-bialgebra
$\Gamma$ in ${}_A\mathrm{BMod}_A(\mathcal{M})$.

The following is one of our main theorems.

\begin{theorem}
[{cf.~Theorem~\ref{thm:right-comodule-main-theorem}}]
\label{thm:main-theorem}
Let $\mathcal{O}^{\otimes}$ be an $\infty$-operad
over a perfect operator category.
Suppose that $\mathcal{M}$ is
an $(\mathcal{O}\times\mathbf{Ass})$-monoidal
$\infty$-category
in which $\mathcal{M}$ has geometric realizations
and the tensor products preserve them
separately in each variable,
and that $A$ is an $(\mathcal{O}\times\mathbf{Ass})$-algebra
object in $\mathcal{M}$.
If $\Gamma$ is an
$(\mathcal{O},\mathbf{Ass})$-bialgebra in 
the mixed $(\mathcal{O},\mathbf{Ass})$-duoidal
$\infty$-category ${}_A\mathrm{BMod}_A(\mathcal{M})$,
then
the $\infty$-category
$\rcomod_{(A,\Gamma)}(\mathcal{M})$
of right $\Gamma$-comodules in
$\mathrm{RMod}_A(\mathcal{M})$ has
the structure of an $\mathcal{O}$-monoidal
$\infty$-category such that
the forgetful functor
$\rcomod_{(A,\Gamma)}(\mathcal{M})
\to \mathrm{RMod}_A(\mathcal{M})$
is strong $\mathcal{O}$-monoidal.
\end{theorem}

The organization of this paper is as follows:
In \S\ref{section:operad}
we review $\infty$-operads over perfect operator categories.
%\if0
%We show that the $\infty$-category of algebra objects
%over the Boardman-Vogt tensor product of
%two $\infty$-operads can be described as
%the $\infty$-category of algebras
%over one $\infty$-operad in the $\infty$-category
%of algebras over the other $\infty$-operad.
%\fi
In \S\ref{section:duoidal-bialgebra}
we recall the definition of
duoidal $\infty$-category.
Mixed duoidal $\infty$-categories
are one of the formulations
of duoidal $\infty$-categories
and they are useful 
when we deal with bialgebras of them.
%We also study bialgebras in a
%mixed duoidal $\infty$-category.
%\if0
%For a mixed duoidal $\infty$-category,
%we can consider bialgebras in it.
%We study the $\infty$-category of bialgebras
%in a mixed duoidal $\infty$-category.
%We show that the $\infty$-category
%nof bialgebras is equivalent
%to the $\infty$-category of monoids in
%the $\infty$-category of comonoids.
%We also show that it is equivalent to
%the $\infty$-category of coalgebras in
%the $\infty$-category of algebras.
%\fi
In \S\ref{section:comodules}
we study multiplicative structures
on bicomodules over bialgebras.
%\if0
%First,
%we recall an approximation 
%of the nonsymmetric $\infty$-operad of bimodules.
%Next,
First,
we set up a generalized bimodule object,
which consists of $\infty$-categories of bimodules,
in the $\infty$-category of $\mathcal{O}$-monoidal
$\infty$-categories and lax $\mathcal{O}$-monoidal functors.
%by looping construction.
From this generalized bimodule object,
we construct an $\infty$-category of bicomodules
which has an $\mathcal{O}$-monoidal structure.
%\fi
In \S\ref{section:right-comodules}
we investigate multiplicative structures
on right comodules over bialgebras.
%\if0
%First,
%we study an approximation 
%of the nonsymmetric $\infty$-operad of right modules.
%Next,
As in the case of bicomodules,
we set up a generalized right module object,
which consists of $\infty$-categories of right modules
and bimodules,
in the $\infty$-category of $\mathcal{O}$-monoidal
$\infty$-categories and lax $\mathcal{O}$-monoidal functors.
We prove Theorem~\ref{thm:main-theorem}
by constructing an $\infty$-category of right comodules
from this generalized right module object.
%\fi

%\newpage
%\input{notation}
%\section{Notation}

\begin{notation}\rm
  
We denote by ${\mathbf\Delta}$ the simplicial indexing category.
We identify a morphism $\sigma: [m]\to [n]$ in
${\mathbf\Delta}$
with a sequence of integers
$(i_0,\ldots,i_m)$ where $0\le i_0\le \cdots\le i_m\le n$.
We abbreviate the sequence $(i,\ldots,i)$ to $(i^n)$,
where $n$ is the number of occurrences of $i$.

Throughout this paper, 
we fix universes $\mathcal{U}\in\mathcal{V}$.
We say that an element of $\mathcal{U}$ is small,
and an element of $\mathcal{V}$ is large.
We denote by
$\cat$ the large $\infty$-category of
small $\infty$-categories,
and 
by $\mathcal{S}$
the large $\infty$-category of small 
$\infty$-groupoids.

For an $\infty$-operad $\mathcal{O}^{\otimes}$
over a perfect operator category,
we denote by $\mathrm{Mon}_{\mathcal{O}}$
the $\infty$-category of small $\mathcal{O}$-monoidal
$\infty$-categories and strong
$\mathcal{O}$-monoidal functors.
We also denote by
$\mathrm{Mon}_{\mathcal{O}}^{\rm (op)lax}$
the $\infty$-category of small $\mathcal{O}$-monoidal
$\infty$-categories and
(op)lax $\mathcal{O}$-monoidal functors.

\if0
Throughout this paper, 
we fix universes $\mathcal{U}\in\mathcal{V}\in\mathcal{W}$.
We say an element of $\mathcal{U}$ is small,
and an element of $\mathcal{V}$ is large,
and an element of $\mathcal{W}$ is very large.
We denote by
$\cat$ the large $\infty$-category of
small $\infty$-categories,
and
by $\wcat$ the very large $\infty$-category
of large $\infty$-categories.
We denote by $\prl$
the subcategory of $\wcat$ spanned
by presentable $\infty$-categories
and colimit-preserving functors.
We denote by $\mathcal{S}$
the large $\infty$-category of small spaces
($\infty$-groupoids).
\fi

\end{notation}

%\newpage
%\input{operad}
\section{$\infty$-operads}\label{section:operad}

In this paper
we use $\infty$-operads over perfect operator
categories in the sense of \cite{Barwick}.
%and their generalizations by \cite{CH}.
The theory of $\infty$-operads over perfect
operator categories is a generalization
of the theory of (symmetric) $\infty$-operads
due to Lurie (\cite[Chapter~2]{Lurie2}).
In this section
we review $\infty$-operads and
generalized $\infty$-operads
over perfect operator categories
and study their properties.

\subsection{$\infty$-operads
over perfect operator categories}

In this subsection
we recall the notions of $\infty$-operads and
generalized $\infty$-operads
over perfect operator categories.

First,
we briefly recall the notion of algebraic pattern
in the sense of \cite[Definition~2.1]{CH}.
Many homotopy-coherent algebraic structures
can be described by algebraic patterns.
An algebraic pattern 
is an $\infty$-category $\mathcal{AP}$
equipped with
\begin{itemize}
\item
an (inert, active) factorization system
$(\mathcal{AP}^{\rm int},\mathcal{AP}^{\rm act})$,
and
\item
a full subcategory $\mathcal{AP}^{\rm el}\subset
\mathcal{AP}^{\rm int}$ of elementary objects.  
\end{itemize}

We can associate an algebraic pattern
to a perfect operator category.
Let $\Phi$ be a perfect operator category
in the sense of \cite[Definition~4.6]{Barwick}.
There is a canonical monad on $\Phi$ 
and we have the Leinster category $\Lambda(\Phi)$
associated to $\Phi$,
which is the Kleisli category of the monad.
The category $\Lambda(\Phi)$ has
an (inert, active) factorization system,
and
we can consider two algebraic patterns
$\Lambda(\Phi)^{\flat}$
and $\Lambda(\Phi)^{\natural}$
by \cite[Example~3.6]{CH}.
The full subcategory $\Lambda(\Phi)^{\flat,\rm el}$
consists only of the final object $*\in\Phi$.
On the other hand,
the full subcategory $\Lambda(\Phi)^{\natural,\rm el}$
contains all objects $E$ such that
there is an inert morphism
$*\to E$ in $\Lambda(\Phi)$.

For an algebraic pattern $\mathcal{AP}$,
we can consider weak $\mathcal{AP}$-Segal fibrations
which are defined in \cite[Definition~9.6]{CH}
in a similar manner to
(symmetric) $\infty$-operads due to Lurie
(\cite[Definition~2.1.1.10]{Lurie2}).
For a perfect operator category $\Phi$,
we define an $\infty$-operad over $\Phi$
to be a weak $\Lambda(\Phi)^{\flat}$-Segal fibration,
and a generalized $\infty$-operad
to be a weak $\Lambda(\Phi)^{\natural}$-Segal fibration.
%in the sense of \cite[Definition~9.6]{CH}.
%(cf.~\cite[Remark~9.8]{CH}).
By \cite[Lemma~9.10]{CH},
$\infty$-operads and generalized $\infty$-operads
over $\Phi$
inherit the structure of algebraic patterns
from $\Lambda(\Phi)^{\flat}$ and $\Lambda(\Phi)^{\natural}$,
respectively.
We denote by $\mathrm{Op}_{\Lambda(\Phi)}$
the $\infty$-category of $\infty$-operads,
and by $\mathrm{Op}_{\Lambda(\Phi)}^{\rm gen}$
the $\infty$-category of generalized $\infty$-operads
over $\Phi$.

For a (generalized) $\infty$-operad
$\mathcal{O}^{\otimes}$ over $\Phi$,
we denote by $\mathcal{O}^{\otimes}_I$ the fiber at $I\in\Phi$.
In particular,
we simply write $\mathcal{O}$
for $\mathcal{O}^{\otimes}_{*}$,
where $*\in\Phi$ is the final object.
For $x\in \mathcal{O}^{\otimes}_I$,
we have a set of inert morphisms
$\{x\to x_i|\ i\in |I|\}$
lying over the inert morphisms
$\{\rho_i: I\to \{i\}|\ i\in |I|\}$.
%Let ${\rm Op}_{\Phi}$ be the $\infty$-category
%of $\infty$-operads over a perfect operator category $\Phi$.
%The $\infty$-category ${\rm Op}_{\Phi}$
%is presentable by \cite[Theorem~7.15]{Barwick}.

Now,
we give examples
of $\infty$-operads over perfect operator categories.
The category of finite sets $\mathbf{F}$
is a perfect operator category
(\cite[Example~1.4]{Barwick}) and
its Leinster category $\Lambda(\mathbf{F})$
is the category of finite pointed sets $\mathrm{Fin}_*$
(\cite[Example~6.5]{Barwick}).
We can recover Lurie's theory of
(symmetric) $\infty$-operads
from the theory of $\infty$-operads
over $\mathbf{F}$
(\cite[Example~7.9]{Barwick}).
Also,
generalized $\infty$-operads over $\mathbf{F}$
coincide with generalized $\infty$-operads
in the sense of \cite[Definition~2.3.2.1]{Lurie2}
(cf.~\cite[Examples~9.8(i)]{CH}).

Nonsymmetric $\infty$-operads are also
formulated by using a perfect operator category.
Recall that $\mathbf{\Delta}$ is the simplicial indexing category.
We can regard $\mathbf{\Delta}^{\rm op}$
as the Leinster category of a perfect operator
category $\mathbf{O}$,
%in the sense of \cite[Definition~4.6]{Barwick},
where $\mathbf{O}$ is the category of ordered finite sets
(\cite[Example~4.9.2]{Barwick}).
We can understand that
$\infty$-operads over $\mathbf{O}$
are nonsymmetric $\infty$-operads
in the sense of \cite[\S3.2]{Gepner-Haugseng}.
Also,
generalized $\infty$-operads over $\mathbf{O}$
are generalized nonsymmetric $\infty$-operads
in the sense of \cite[Definition~2.4.1]{Gepner-Haugseng}
(cf.~\cite[Examples~9.8(ii)]{CH}).
%\end{definition}

\subsection{Monoidal $\infty$-categories}

In this subsection
we recall that there are several equivalent formulations of
monoidal $\infty$-categories.
In particular,
we recall interaction between monoidal
$\infty$-categories and 
monoid objects in the $\infty$-category
of small $\infty$-categories.

Let $\Phi$ be a perfect operator category
and let $\mathcal{O}^{\otimes}$
be a (generalized) $\infty$-operad over $\Phi$.
We abbreviate 
the slice category
$(\mathrm{Op}_{\Lambda(\Phi)}^{\rm (gen)})_{/\mathcal{O}^{\otimes}}$
as ${\rm Op}_{\mathcal{O}}^{\rm (gen)}$.
%the $\infty$-category
%${\rm Op}_{\Phi/\mathcal{O}^{\otimes}}$
%of $\infty$-operads
%over $\mathcal{O}^{\otimes}$.
%The $\infty$-category ${\rm Op}_{\mathcal{O}}$
%is presentable by \cite[Proposition~5.5.3.10]{Lurie1}.
%For $\infty$-operads
%$\mathcal{O}^{\otimes}$ and $\mathcal{P}^{\otimes}$
%over perfect operator categories,
For objects $\mathcal{X}^{\otimes}$ and $\mathcal{Y}^{\otimes}$
of $\mathrm{Op}_{\mathcal{O}}^{\rm (gen)}$,
we write 
$\mathrm{Alg}_{\mathcal{X}/\mathcal{O}}(\mathcal{Y})$
for the full subcategory
of ${\rm Fun}_{\mathcal{O}^{\otimes}}
(\mathcal{X}^{\otimes},\mathcal{Y}^{\otimes})$
spanned by morphisms of $\mathrm{Op}_{\mathcal{O}}^{\rm (gen)}$.
When $\mathcal{X}^{\otimes}=\mathcal{O}^{\otimes}$,
we simply write $\mathrm{Alg}_{\mathcal{O}}(\mathcal{Y})$.

Now,
we suppose that $\mathcal{O}^{\otimes}$
is an $\infty$-operad.
An $\mathcal{O}$-monoidal $\infty$-category
is a cocartesian fibration
of $\infty$-operads
$\mathcal{M}^{\otimes}\to\mathcal{O}^{\otimes}$
over $\Phi$.
For simplicity,
we also say that 
$\mathcal{M}$ is an $\mathcal{O}$-monoidal $\infty$-category.
We write $\monlax$
for the full subcategory of ${\rm Op}_{\mathcal{O}}$
spanned by $\mathcal{O}$-monoidal $\infty$-categories.
We say that a morphism in $\monlax$
is a lax $\mathcal{O}$-monoidal functor.

A morphism in $\mathrm{Op}_{\mathcal{O}}$
between $\mathcal{O}$-monoidal $\infty$-categories
is said to be a strong $\mathcal{O}$-monoidal functor
if it preserves cocartesian morphisms.
We write $\mons$
for the wide subcategory of $\monlax$
spanned by strong $\mathcal{O}$-monoidal functors.

Next,
we will formulate oplax $\mathcal{O}$-monoidal functors
between $\mathcal{O}$-monoidal $\infty$-categories.
Let $(-)^{\rm op}: \cat\to\cat$
be the functor which assigns to an $\infty$-category
$\mathcal{C}$ its opposite $\infty$-category $\mathcal{C}^{\rm op}$.
We say that a functor of $\infty$-categories
$\mathcal{C}\to \mathcal{O}^{\otimes,\rm op}$
is an opposite $\infty$-operad over $\mathcal{O}^{\otimes,\rm op}$
if its opposite functor $\mathcal{C}^{\rm op}\to
\mathcal{O}^{\otimes}$ is a map of $\infty$-operads.
For opposite $\infty$-operads $\mathcal{C}$ and $\mathcal{D}$
over $\mathcal{O}^{\otimes, \rm op}$,
a map $\mathcal{C}\to \mathcal{D}$ over $\mathcal{O}^{\otimes,\rm op}$
is a morphism of opposite $\infty$-operads
over $\mathcal{O}^{\otimes,\rm op}$ 
if its opposite $\mathcal{C}^{\rm op}\to\mathcal{D}^{\rm op}$
is a morphism of $\infty$-operads
over $\mathcal{O}^{\otimes}$.
We denote by $\mathrm{Op}_{\mathcal{O}}^{\vee}$
the subcategory of $\mathrm{Cat}_{\infty/\mathcal{O}^{\otimes,\rm op}}$
spanned by opposite $\infty$-operads
over $\mathcal{O}^{\otimes,\rm op}$
and their morphisms.

We say that an opposite $\infty$-operad
$\mathcal{C}\to\mathcal{O}^{\otimes,\rm op}$
over $\mathcal{O}^{\otimes,\rm op}$
is a cartesian fibration of opposite $\infty$-operads
if its opposite 
$\mathcal{C}^{\rm op}\to\mathcal{O}^{\otimes}$
is a cocartesian fibration of $\infty$-operads.
We write $\mons^{\rm oplax}$
for the full subcategory of $\mathrm{Op}_{\mathcal{O}}^{\vee}$
spanned by cartesian fibrations of opposite $\infty$-operads.
We say that a morphism of $\monoplax$
is an oplax $\mathcal{O}$-monoidal functor. 

We notice that the functor
$(-)^{\rm op}$ induces an equivalence of
$\infty$-categories
$(-)^{\rm op}: \monlax\to \monoplax$.
A morphism of cartesian fibrations
of opposite $\infty$-operads over $\mathcal{O}^{\otimes, \rm op}$
is said to be a strong $\mathcal{O}$-monoidal functor
if it preserves cartesian morphisms.
We denote by $\mons^{\vee}$
the wide subcategory of $\mons^{\rm oplax}$
spanned by strong $\mathcal{O}$-monoidal functors.
The equivalence $(-)^{\rm op}$ restricts to
an equivalence of $\infty$-categories
$(-)^{\rm op}: \mons\to\mons^{\vee}$.

Next,
we recall monoid objects
in an $\infty$-category.
Let $\mathcal{C}$ be an $\infty$-category
with finite products.
We say that a functor
$M: \mathcal{O}^{\otimes}\to\mathcal{C}$
is an $\mathcal{O}$-monoid object of $\mathcal{C}$
if the Segal map
\[ M(x)\longrightarrow \prod_{i\in |I|}M(x_i) \]
induced by inert morphisms
$\{x\to x_i|\ i\in |I|\}$
is an equivalence in $\mathcal{C}$
for any $x\in\mathcal{O}^{\otimes}$.
%over $I\in\Phi$.
We define an $\infty$-category
${\rm Mon}_{\mathcal{O}}(\mathcal{C})$
to be the full subcategory
of ${\rm Fun}(\mathcal{O}^{\otimes},\mathcal{C})$
spanned by $\mathcal{O}$-monoid objects.

We notice that
an $\mathcal{O}$-monoid object coincides
with a Segal $\mathcal{O}^{\otimes,\flat}$-object
in the sense of \cite[Definition~2.7]{CH},
where $\mathcal{O}^{\otimes,\flat}$
is the algebraic pattern
coming from $\Lambda(\Phi)^{\flat}$.
When $\mathcal{O}^{\otimes}$ is a generalized
$\infty$-operad
and $\mathcal{C}$ is an $\infty$-category
with finite limits,
we define a generalized $\mathcal{O}$-monoid object
in $\mathcal{C}$
to be a Segal $\mathcal{O}^{\otimes,\natural}$-object,
where $\mathcal{O}^{\otimes,\natural}$
is the algebraic pattern
coming from $\Lambda(\Phi)^{\natural}$.
We denote by $\mathrm{Mon}_{\mathcal{O}}^{\rm gen}(\mathcal{C})$
the full $\infty$-category of
$\mathrm{Fun}(\mathcal{O}^{\otimes},\mathcal{C})$
spanned by generalized $\mathcal{O}$-monoid objects.

For example,
we set $\mathbf{Ass}^{\otimes}=\mathbf{\Delta}^{\rm op}$,
and regard the identity map
$\mathbf{Ass}^{\otimes}\to\mathbf{\Delta}^{\rm op}$
as an $\infty$-operad over $\mathbf{O}$.
An $\mathbf{Ass}$-monoid
object of $\mathcal{C}$
is a functor $M:
{\mathbf\Delta}^{\rm op}\to \mathcal{C}$
%\mathbf{Ass}^{\otimes}\to \mathcal{C}$
such that the map
$M([n])\to M([1])^{\times n}$
induced by inert morphisms $[1]\to [n]$ in ${\mathbf\Delta}$
is an equivalence
for all $n\ge 0$.
An $\mathbf{Ass}$-monoid object
is simply called a monoid object.

We can also regard the identity map
$\mathbf{Ass}^{\otimes}\to{\mathbf\Delta}^{\rm op}$
as a generalized $\infty$-operad.
A generalized $\mathbf{Ass}$-monoid
object is a functor
$M: {\mathbf\Delta}^{\rm op}\to\mathcal{C}$
such that the map
$M([n])\to
   M([1])^{\times_{M([0])}n}$
%$M([n])\to
%   M([1])\times_{M([0])}\cdots\times_{M([0])}M([1])$
induced by inert morphisms $[i]\to [n]$
for $i=0,1$ in $\mathbf{\Delta}$
is an equivalence 
for all $n\ge 0$.
A generalized $\mathbf{Ass}$-monoid object
is also called a category object.

Now,
we consider a relationship
between algebra objects and monoid objects.
When $\mathcal{C}$
is an $\infty$-category 
with finite products,
we can promote $\mathcal{C}$ to a cartesian symmetric
monoidal $\infty$-category
$\mathcal{C}^{\times}\to \mathrm{Fin}_*$
(cf.~\cite[\S2.4.1]{Lurie2}).
There is an operator morphism
$u: \Phi\to \mathbf{F}$
given by $I\mapsto |I|$
(\cite[Example~1.11]{Barwick}),
which induces an adjunction
of $\infty$-categories
\[ u_!: \mathrm{Op}_{\Lambda(\Phi)}\rightleftarrows
        \mathrm{Op}_{\mathrm{Fin}_*}: u^* \]
by \cite[Proposition~7.18]{Barwick}.
In the same way as
the proof of \cite[Proposition~2.4.17]{Lurie2},
we can show that 
there are equivalences of $\infty$-categories
%\begin{equation}\label{eq:Mon-Alg-equivalence}
\[
\mathrm{Mon}_{\mathcal{O}}(\mathcal{C})\simeq
\mathrm{Alg}_{\mathcal{O}}(u^*\mathcal{C}^{\times})
\simeq
\mathrm{Alg}_{u_!\mathcal{O}}(\mathcal{C}^{\times}).
\]
%\end{equation}

Finally,
we recall an interaction of monoidal $\infty$-categories
and monoid objects in the $\infty$-category
of small $\infty$-categories $\cat$.
Since
%the $\infty$-category of small $\infty$-categories
$\cat$ has finite products,
we can consider an $\infty$-category
of $\mathcal{O}$-monoid objects ${\rm Mon}_{\mathcal{O}}(\cat)$.
To an $\mathcal{O}$-monoid object
$X: \mathcal{O}^{\otimes}\to\cat$,
we can associate a cocartesian fibration
of $\infty$-operads $X^{\otimes}\to\mathcal{O}^{\otimes}$
by unstraightening.
This correspondence induces an equivalence
of $\infty$-categories between 
$\mathrm{Mon}_{\mathcal{O}}(\cat)$
and $\mathrm{Mon}_{\mathcal{O}}$.
%$\mathrm{Alg}_{\mathcal{O}}(u^*\cat^{\times})$.
%by (\ref{eq:Mon-Alg-equivalence}).

%\newpage

\subsection{Sections of two-variable fibrations}

In this subsection
we study $\infty$-categories of sections
of two-variable fibrations
for later use.

Let $p=(f,g): X\to A\times B$
be a categorical fibration of $\infty$-categories,
%and let $B'\in{\rm Cat}_{\infty/B}$.
and let $\psi: B'\to B$ be a map of $\infty$-categories.
We define a map of $\infty$-categories
$\xi: F_{B'/B}(X)\to A$
by the following pullback diagram
\[ \xymatrix{
    F_{B'/B}(X)\ar[r]\ar[d]_{\xi}&
    {\rm Fun}(B',X)\ar[d]^{(f,g)_*} \\
    A\ar[r] & 
    {\rm Fun}(B',A\times B),\\
}\]
where the bottom horizontal arrow
is an adjoint of the map
%$A\times B'\to A\times B$.
${\rm id}_A\times\psi: A\times B'\to A\times B$.
We notice that there is a natural equivalence
\[ {\rm Map}_{{\rm Cat}_{\infty/A}}(K,F_{B'/B}(X))\simeq
   {\rm Map}_{{\rm Cat}_{\infty/A\times B}}
   (K\times B',X) \]
of $\infty$-groupoids
for any $K\in {\rm Cat}_{\infty/A}$.
In particular,
an object of $F_{B'/B}(X)$ is identified
with a map $B'\to X_a$ of $\infty$-categories
over $B$ for some $a\in A$,
where $X_a$ is the fiber of $f$ at $a$.

We will show that
$\xi$ is a cocartesian fibration under some conditions.
%We denote by $\pi_A$ the projection $A\times B\to A$.

\begin{lemma}\label{lemma:xi-cocartesian-fibration}
If $f: X\to A$ is a cocartesian fibration
and the map $g$ takes
$f$-cocartesian edges to equivalences,
then $\xi: F_{B'/B}(X)\to A$
is a cocartesian fibration. 
\end{lemma}

\proof
The pair of maps $(f,g)$ induces
a categorical fibration
$(f_*,g_*):
{\rm Fun}(B',X)\to {\rm Fun}(B',A)\times {\rm Fun}(B',B)$,
where $f_*: {\rm Fun}(B',X)\to {\rm Fun}(B',A)$
is a cocartesian fibration and 
$g_*$ takes $f_*$-cocartesian edges
to equivalences.
By taking the fiber at $\psi$,
we obtain a cocartesian fibration
${\rm Fun}(B',X)_{\psi}\to {\rm Fun}(B',A)$,
where ${\rm Fun}(B',X)_{\psi}=
{\rm Fun}(B',X)\times_{{\rm Fun}(B',B)}\{\psi\}$.
The lemma follows from the fact that
we have a pullback diagram
\[ \xymatrix{
    F_B(B',X)\ar[r]\ar[d] & 
    {\rm Fun}(B',X)_{\psi}\ar[d]\\
    A\ar[r] &
    {\rm Fun}(B',A),\\
   }\]
where the bottom horizontal arrow is an adjoint
of the projection $A\times B'\to A$.
\qed

\begin{remark}\rm
Under the assumption of
Lemma~\ref{lemma:xi-cocartesian-fibration},
the map $p=(f,g): X\to A\times B$
induces a functor
$A\to {\rm Cat}_{\infty/B}$
by straightening.
In this situation
the cocartesian fibration
$\xi: F_B(B', X)\to A$
is associated to the composite functor
\[ A\longrightarrow {\rm Cat}_{\infty/B}
   {\longrightarrow}\cat \] 
by unstraightening,
where the second arrow is
${\rm Fun}_B(B',- )$.
\end{remark}

\begin{remark}\label{remark:xi-cartesian-fibration}\rm
Suppose that $f: X\to A$ is a cartesian fibration
and that $g$ takes $f$-cartesian edges
to equivalences.
By a dual argument
to the proof of Lemma~\ref{lemma:xi-cocartesian-fibration},
we see that
$\xi: F_{B'/B}(X)\to A$
is a cartesian fibration.
\end{remark}

\subsection{Algebra objects over
products of $\infty$-operads}
%Boardman-Vogt tensor products}

In this subsection
we study $\infty$-categories
of algebra objects over products of $\infty$-operads.

Let $\Phi$ and $\Psi$ be perfect operator categories.
The product $\Phi\times\Psi$ is also a perfect operator category,
and its Leinster category $\Lambda(\Phi\times\Psi)$
is isomorphic to the product $\Lambda(\Phi)\times\Lambda(\Psi)$.
Let 
$\mathcal{O}^{\otimes}$ and $\mathcal{P}^{\otimes}$
be $\infty$-operads over $\Phi$ and $\Psi$, respectively.
The product
$\mathcal{O}^{\otimes}\times\mathcal{P}^{\otimes}$
is an $\infty$-operad over $\Phi\times\Psi$.

First,
we study monoid objects over
$\mathcal{O}^{\otimes}\times\mathcal{P}^{\otimes}$.
For an $\infty$-category $\mathcal{C}$
with finite products,
the $\infty$-category ${\rm Mon}_{\mathcal{P}}(\mathcal{C})$
of $\mathcal{P}$-monoid objects of $\mathcal{C}$
also has finite products.
Thus,
we can consider an $\mathcal{O}$-monoid object
of ${\rm Mon}_{\mathcal{P}}(\mathcal{C})$.
%We have the following lemmas.
%By \cite[Corollary~9.10.1]{Barwick},
%we obtain the following lemma.

\begin{lemma}[{cf.~\cite[Example~5.7]{CH}}]
\label{lemma:monoidal-category-over-tensor-product}
Let $\mathcal{C}$ be an $\infty$-category with finite products.
For $\infty$-operads
$\mathcal{O}^{\otimes}$ and $\mathcal{P}^{\otimes}$
over perfect operator categories,
we have equivalences of $\infty$-categories
\[ {\rm Mon}_{\mathcal{O}}({\rm Mon}_{\mathcal{P}}(\mathcal{C}))
   \simeq
   {\rm Mon}_{\mathcal{O}\times\mathcal{P}}(\mathcal{C})
   \simeq
   {\rm Mon}_{\mathcal{P}}({\rm Mon}_{\mathcal{O}}(\mathcal{C})). \]
\end{lemma}

\proof
This follows by observing that
${\rm Mon}_{\mathcal{O}\times\mathcal{P}}(\mathcal{C})$
is a full subcategory
of ${\rm Fun}(\mathcal{O}^{\otimes}\times\mathcal{P}^{\otimes},
\mathcal{C})$
spanned by those functors $M$
such that the induced maps
$M(x,y)\to \prod_{i\in |I|}M(x_i,y)$
and $M(x,y)\to \prod_{j\in |J|}M(x,y_j)$
are equivalences in $\mathcal{C}$
for any $(x,y)\in\mathcal{O}^{\otimes}
\times\mathcal{P}^{\otimes}$.
%where $\{x\to x_i|\ {i\in |I|}\}$
%and $\{y\to y_j|\ j\in |J|\}$
%are sets of inert morphisms lying over
%inert morphisms $\{I\to \{i\}|\ i\in |I|\}$
%and $\{J\to \{j\}|\ i\in |J|\}$,
%respectively.
\qed

%@\bigskip

Next,
we study $\infty$-category of algebra objects
in monoidal $\infty$-categories over
products of $\infty$-operads.
We suppose that
$\mathcal{A}^{\otimes}\to\mathcal{O}^{\otimes}$
and
$\mathcal{B}^{\otimes}\to\mathcal{P}^{\otimes}$
are maps of $\infty$-operads
over $\Phi$ and $\Psi$,
respectively.
   
\begin{lemma}\label{lemma:algebra-monoidal-structure}
For an $(\mathcal{O}\times\mathcal{P})$-monoidal
$\infty$-category $\mathcal{C}$,
the $\infty$-category
${\rm Alg}_{\mathcal{B}/\mathcal{P}}(\mathcal{C})$
%${\rm Alg}_{\mathcal{B}/\mathcal{P}}^{\rm Mon}(\mathcal{C})$
has the structure of an $\mathcal{O}$-monoidal $\infty$-category.
\end{lemma}

\proof
By Lemma~\ref{lemma:monoidal-category-over-tensor-product},
the $(\mathcal{O}\times\mathcal{P})$-monoidal $\infty$-category
$\mathcal{C}$ determines a functor
$\overline{\mathcal{C}}:
\mathcal{O}^{\otimes}\to {\rm Mon}_{\mathcal{P}}(\cat)$
that is an $\mathcal{O}$-monoid object.
Since the functor
${\rm Alg}_{\mathcal{B}/\mathcal{P}}(-)$
%${\rm Alg}_{\mathcal{B}/\mathcal{P}}^{\rm Mon}(-)$
preserves finite products,
the composite
of $\overline{\mathcal{C}}$ with
${\rm Alg}_{\mathcal{B}/\mathcal{P}}(-)$
%${\rm Alg}_{\mathcal{B}/\mathcal{P}}^{\rm Mon}(-)$
is an $\mathcal{O}$-monoid object.
This gives us an $\mathcal{O}$-monoidal structure on
${\rm Alg}_{\mathcal{B}/\mathcal{P}}(\mathcal{C})$.
%${\rm Alg}_{\mathcal{B}/\mathcal{P}}^{\rm Mon}(\mathcal{C})$.
\qed

%@\bigskip

We will show that an analogue
of Lemma~\ref{lemma:monoidal-category-over-tensor-product}
holds for
${\rm Alg}_{(\mathcal{A}\times\mathcal{B})/
(\mathcal{O}\times\mathcal{P})}(\mathcal{C})$.

\begin{proposition}\label{prop:double-algebra-decomposition}
For an $(\mathcal{O}\times\mathcal{P})$-monoidal
$\infty$-category $\mathcal{C}$,
there are equivalences of $\infty$-categories
\[  {\rm Alg}_{\mathcal{A}/\mathcal{O}}
    ({\rm Alg}_{\mathcal{\mathcal{B}/\mathcal{P}}}(\mathcal{C}))
    \simeq
    {\rm Alg}_{(\mathcal{A}\times\mathcal{B})/
    (\mathcal{O}\times\mathcal{P})}(\mathcal{C})
   \simeq
       {\rm Alg}_{\mathcal{B}/\mathcal{P}}
       ({\rm Alg}_{\mathcal{A}/\mathcal{O}}(\mathcal{C})).\]
\end{proposition}

\proof
We will prove the first equivalence.
The second equivalence can be proved
in a similar manner.
%First,
%we notice that 
%for any $\mathcal{R}$-monoidal $\infty$-category $\mathcal{D}$,
%we have a fully faithful functor
%${\rm Alg}_{\mathcal{Q}/\mathcal{R}}(\mathcal{D})
%   \to
%   {\rm Fun}_{\mathcal{R}^{\otimes}}(\mathcal{Q}^{\otimes},
%   \mathcal{D}^{\otimes})$.
By Lemma~\ref{lemma:algebra-monoidal-structure},
$\mathrm{Alg}_{\mathcal{B}/\mathcal{P}}(\mathcal{C})$
%$\mathrm{Alg}_{\mathcal{B}/\mathcal{P}}^{\rm Mon}(\mathcal{C})$
has the structure of an $\mathcal{O}$-monoidal $\infty$-category.
We denote by
$\mathrm{Alg}_{\mathcal{B}/\mathcal{P}}(\mathcal{C})^{\otimes}
\to\mathcal{O}^{\otimes}$
%$\mathrm{Alg}_{\mathcal{B}/\mathcal{P}}^{\rm Mon}(\mathcal{C})^{\otimes}
%\to\mathcal{O}^{\otimes}$
the corresponding cocartesian fibration
of $\infty$-operads
by unstraightening.
By Lemma~\ref{lemma:xi-cocartesian-fibration},
we have a cocartesian fibration
$F_{\mathcal{B}^{\otimes}/\mathcal{P}^{\otimes}}
(\mathcal{C}^{\otimes})\to \mathcal{O}^{\otimes}$.
We can verify that it is a map of $\infty$-operads.
The evaluation map
$\mathrm{Alg}_{\mathcal{B}/\mathcal{P}}
(\mathcal{C})^{\otimes}
\times\mathcal{B}^{\otimes}\to\mathcal{C}^{\otimes}$
%$\mathrm{Alg}_{\mathcal{B}/\mathcal{P}}^{\rm Mon}
%(\mathcal{C})^{\otimes}
%\times\mathcal{B}^{\otimes}\to\mathcal{C}^{\otimes}$
induces a map of $\infty$-operads
$\mathrm{Alg}_{\mathcal{B}/\mathcal{P}}
(\mathcal{C})^{\otimes}
\to F_{\mathcal{B}/\mathcal{P}}(\mathcal{C}^{\otimes})$
%$\mathrm{Alg}_{\mathcal{B}/\mathcal{P}}^{\rm Mon}
%(\mathcal{C})^{\otimes}
%\to F_{\mathcal{B}/\mathcal{P}}(\mathcal{C}^{\otimes})$
over $\mathcal{O}^{\otimes}$.
We notice that the induced map on fibers
is identified with the fully faithful functor
%$\mathrm{Al}_{\mathcal{B}/\mathcal{P}}^{\rm Mon}
$\mathrm{Al}_{\mathcal{B}/\mathcal{P}}
(\mathcal{C}_I)
\to \mathrm{Fun}_{\mathcal{P}^{\otimes}}
(\mathcal{B}^{\otimes},\mathcal{C}^{\otimes}_I)$
for any $I\in\mathcal{O}^{\otimes}$.
This implies that the composite 
${\rm Alg}_{\mathcal{A}/\mathcal{O}}
 ({\rm Alg}_{\mathcal{B}/\mathcal{P}}(\mathcal{C}))
   \to
   {\rm Alg}_{\mathcal{A}/\mathcal{O}}
   (F_{\mathcal{B}^{\otimes}/\mathcal{P}^{\otimes}}
   (\mathcal{C}^{\otimes}))
\to
   {\rm Fun}_{\mathcal{O}^{\otimes}}
   (\mathcal{A}^{\otimes},
   F_{\mathcal{B}^{\otimes}/\mathcal{P}^{\otimes}}
   (\mathcal{C}^{\otimes}))$
\if0
${\rm Alg}_{\mathcal{A}/\mathcal{O}}^{\rm Mon}
   ({\rm Alg}_{\mathcal{B}/\mathcal{P}}^{\rm Mon}(\mathcal{C}))
   \to
   {\rm Alg}_{\mathcal{A}/\mathcal{O}}^{\rm Mon}
   (F_{\mathcal{B}^{\otimes}/\mathcal{P}^{\otimes}}
   (\mathcal{C}^{\otimes}))
\to
   {\rm Fun}_{\mathcal{O}^{\otimes}}
   (\mathcal{A}^{\otimes},
   F_{\mathcal{B}^{\otimes}/\mathcal{P}^{\otimes}}
   (\mathcal{C}^{\otimes}))$
\fi
is fully faithful,
where the right hand side is equivalent to
${\rm Fun}_{\mathcal{O}^{\otimes}\times\mathcal{P}^{\otimes}}
   (\mathcal{A}^{\otimes}\times\mathcal{B}^{\otimes},
    \mathcal{C}^{\otimes})$.
We also have a fully faithful functor 
${\rm Alg}_{(\mathcal{A}\times\mathcal{B})/
  (\mathcal{O}\times\mathcal{P})}(\mathcal{C})
  \to
  {\rm Fun}_{\mathcal{O}^{\otimes}\times\mathcal{P}^{\otimes}}
  (\mathcal{A}^{\otimes}\times\mathcal{B}^{\otimes},
  \mathcal{C}^{\otimes})$.
\if0
${\rm Alg}_{(\mathcal{A}\times\mathcal{B})/
  (\mathcal{O}\times\mathcal{P})}^{\rm Mon}(\mathcal{C})
  \to
  {\rm Fun}_{\mathcal{O}^{\otimes}\times\mathcal{P}^{\otimes}}
  (\mathcal{A}^{\otimes}\times\mathcal{B}^{\otimes},
  \mathcal{C}^{\otimes})$.
\fi
The proposition follows
by observing that the essential images of
${\rm Alg}_{(\mathcal{A}\times\mathcal{B})/
     (\mathcal{O}\times\mathcal{P})}(\mathcal{C})$
%${\rm Alg}_{(\mathcal{A}\times\mathcal{B})/
%     (\mathcal{O}\times\mathcal{P})}^{\rm Mon}(\mathcal{C})$
and
${\rm Alg}_{\mathcal{A}/\mathcal{O}}
   ({\rm Alg}_{\mathcal{B}/\mathcal{P}}(\mathcal{C}))$
%${\rm Alg}_{\mathcal{A}/\mathcal{O}}^{\rm Mon}
%   ({\rm Alg}_{\mathcal{B}/\mathcal{P}}^{\rm Mon}(\mathcal{C}))$
coincide in
$ {\rm Fun}_{\mathcal{O}^{\otimes}\times\mathcal{P}^{\otimes}}
  (\mathcal{A}^{\otimes}\times\mathcal{B}^{\otimes},
  \mathcal{C}^{\otimes})$.
\qed

%\newpage
%\input{change-operad}

%\newpage
%\input{virtual-duoidal}
\section{Duoidal $\infty$-categories and bialgebras}
\label{section:duoidal-bialgebra}

In this section
we recall the definition of duoidal $\infty$-categories
and study bialgebras in duoidal $\infty$-categories.
There are three formulations
of duoidal $\infty$-categories
according to which kinds of functors we take.
In order to formulate bialgebras
in duoidal $\infty$-categories
it is convenient to use mixed duoidal $\infty$-categories.
%Next,
%we study $\infty$-categories of sections of
%mixed duoidal $\infty$-categories
%since bialgebras are defined to be sections
%satisfying some conditions.
We will introduce an $\infty$-categories of bialgebras
in a mixed duoidal $\infty$-category
and show that
it can be described by $\infty$-categories
of algebras and coalgebras.

\subsection{Mixed duoidal $\infty$-categories}

In this subsection
we recall the notion of mixed duoidal $\infty$-categories.
Mixed duoidal $\infty$-categories
are one of the formulations of duoidal $\infty$-categories
which are adapted for bilax monoidal functors.
Mixed duoidal $\infty$-categories
are defined by using mixed fibrations.

%First,
%we recall the definition of mixed fibrations.

\begin{definition}[{cf.~\cite[Definition~4.15]{Torii1}}]\rm
Let $p: X\to A\times B$ be a categorical
fibration of $\infty$-categories.
We say that $p$ is a mixed fibration
if it satisfies the following conditions:

\begin{enumerate}
\item
The composite $\pi_A\circ p: X\to A$
is a cocartesian fibration,
and the map $p$ takes $(\pi_A\circ p)$-cocartesian
edges to $\pi_A$-cocartesian edges,
where $\pi_A: A\times B\to A$ is the projection.

\item
The composite $\pi_B\circ p: X\to B$
is a cartesian fibration,
and the map $p$ takes $(\pi_B\circ p)$-cartesian
edges to $\pi_B$-cartesian edges,
where $\pi_B: A\times B\to B$ is the projection.
  
\end{enumerate}
\end{definition}

\begin{remark}\rm
Mixed fibrations are said to be
curved orthofibrations in \cite{HHLN1}
and
lax two-sided fibrations in \cite{Stefanich}.
\end{remark}

Let $\mathcal{P}^{\otimes}$ and
$\mathcal{Q}^{\boxtimes}$ be $\infty$-operads
over perfect operator categories
$\Phi$ and $\Psi$,
%in the sense of \cite{Barwick},
respectively.
Suppose that we have a mixed fibration 
\[ p: X\longrightarrow
   \mathcal{P}^{\otimes}\times \mathcal{Q}^{\boxtimes,\rm op}.\]

\begin{definition}[{cf.~\cite[Definition~5.10]{Torii1}
and \cite[Definition~3.11]{Torii2}}]
\rm
We say that 
$p$ is a mixed $(\mathcal{P},\mathcal{Q})$-duoidal
$\infty$-category if
it satisfies the following conditions:

\begin{enumerate}

\item
For any $x\in\mathcal{P}^{\otimes}$,
the Segal morphism
$X_x\to \prod_{i\in |I|}^{\mathcal{Q}^{\boxtimes,\rm op}}X_{x_i}$
is an equivalence in ${\rm Cat}_{\infty/\mathcal{Q}^{\boxtimes,\rm op}}$.

\item
For any $y\in\mathcal{Q}^{\boxtimes, \rm op}$,
the Segal morphism
$X_y\to \prod_{j\in |J|}^{\mathcal{P}^{\otimes}}X_{y_j}$
is an equivalence in ${\rm Cat}_{\infty/\mathcal{P}^{\otimes}}$.
  
\end{enumerate}
\end{definition}

\begin{definition}[{cf.~\cite[Definition~5.13]{Torii1}
and \cite[Definition~3.11]{Torii2}}]
\rm
Let $p: X\to \mathcal{P}^{\otimes}\times
\mathcal{Q}^{\boxtimes,\rm op}$ and
$q: Y\to \mathcal{P}^{\otimes}\times
\mathcal{Q}^{\boxtimes,\rm op}$ be
mixed $(\mathcal{P},\mathcal{Q})$-duoidal
$\infty$-categories.
A functor $f: X\to Y$ over $\mathcal{P}^{\otimes}\times
\mathcal{Q}^{\boxtimes,\rm op}$
is said to be a bilax $(\mathcal{P},\mathcal{Q})$-monoidal
functor if it preserves cocartesian edges
over inert morphisms of $\mathcal{P}^{\otimes}$
and cartesian edges over inert morphisms
of $\mathcal{Q}^{\boxtimes,\rm op}$.
We denote by
\[ \mathrm{Duo}_{(\mathcal{P},\mathcal{Q})}^{\rm bilax} \]
the subcategory of
$\mathrm{Cat}_{\infty/\mathcal{P}^{\otimes}\times\mathcal{Q}^{\boxtimes,\rm op}}$
spanned by mixed $(\mathcal{P},\mathcal{Q})$-duoidal
$\infty$-categories and
bilax $(\mathcal{P},\mathcal{Q})$-monoidal functors.
\end{definition}

\begin{remark}\label{remark:equivalence-mixed-duoidal-bilax}
\rm
By \cite[Theorem~5.17 and  Remark~5.18]{Torii1}
and \cite[Theorem~3.12]{Torii2},
there are equivalences of $\infty$-categories
\[ \mathrm{Mon}_{\mathcal{P}}^{\rm lax}
   (\mathrm{Mon}_{\mathcal{Q}}^{\rm oplax})
   \simeq
   \mathrm{Duo}_{(\mathcal{P},\mathcal{Q})}^{\rm bilax}
   \simeq
   \mathrm{Mon}_{\mathcal{Q}}^{\rm oplax}
   (\mathrm{Mon}_{\mathcal{P}}^{\rm lax}).\]
In particular,
a (mixed) $(\mathcal{P},\mathcal{Q})$-duoidal
$\infty$-category can be identified with
a $\mathcal{P}$-monoid object
of $\mathrm{Mon}_{\mathcal{Q}}^{\rm oplax}$,
and a $\mathcal{Q}$-monoid object
of $\mathrm{Mon}_{\mathcal{P}}^{\rm lax}$.
\end{remark}

In a similar manner to
\cite[Proposition~5.12 and Remark~5.16]{Torii1},
we have the following proposition.

\begin{proposition}\label{prop:duality-mixed-duoidal}
The functor $(-)^{\rm op}:\cat\to\cat$
induces an equivalence of $\infty$-categories
\[ (-)^{\rm op}:
   \mathrm{Duo}_{(\mathcal{P},\mathcal{Q})}^{\rm bilax}
   \stackrel{\simeq}{\longrightarrow}
   \mathrm{Duo}_{(\mathcal{Q},\mathcal{P})}^{\rm bilax}.\]
\end{proposition}

%\newpage

\subsection{Bialgebras in mixed duoidal $\infty$-categories}

In this subsection
we will introduce an $\infty$-category
of bialgebras in a mixed duoidal
$\infty$-category.
We will show that 
algebras and coalgebras in a mixed duoidal $\infty$-category
form monoidal $\infty$-categories, respectively,
and that 
the $\infty$-category of bialgebras
can be described by
$\infty$-categories of algebras and coalgebras.

In this subsection
we suppose that we have a mixed
$(\mathcal{P},\mathcal{Q})$-duoidal
$\infty$-category
$p: X\to
\mathcal{P}^{\otimes}\times \mathcal{Q}^{\boxtimes,\rm op}$.
Furthermore,
we suppose that
$\mathcal{A}^{\otimes}\to\mathcal{P}^{\otimes}$
is a map of $\infty$-operads over $\Phi$,
and that
$\mathcal{B}^{\boxtimes}\to\mathcal{Q}^{\boxtimes}$
is a map of $\infty$-operads over $\Psi$.

First,
we will construct an $\infty$-category
${\rm Alg}_{\mathcal{A}/\mathcal{P}}(X)^{\boxtimes,\vee}$
over $\mathcal{Q}^{\boxtimes,\rm op}$
which consists of $\mathcal{A}$-algebra objects
in $X_y$
for some $y\in\mathcal{Q}^{\boxtimes, \rm op}$.
Next,
we will show that
${\rm Alg}_{\mathcal{A}/\mathcal{P}}(X)^{\boxtimes, \vee}
\to \mathcal{Q}^{\boxtimes,\rm op}$
is a cartesian fibration of opposite
$\infty$-operads.

By Remark~\ref{remark:xi-cartesian-fibration},
we have a cartesian fibration
$F_{\mathcal{A}^{\otimes}/\mathcal{P}^{\otimes}}(X)
\to \mathcal{Q}^{\boxtimes, \rm op}$,
where an object of
$F_{\mathcal{A}^{\otimes}/\mathcal{P}^{\otimes}}(X)$
is identified with a map
$\mathcal{A}^{\otimes}\to X_y$
over $\mathcal{P}^{\otimes}$
for some $y\in\mathcal{Q}^{\boxtimes,\rm op}$.

\begin{definition}\rm
We define an $\infty$-category
\[ {\rm Alg}_{\mathcal{A}/\mathcal{P}}(X)^{\boxtimes, \vee} \]
to be the full subcategory of
$F_{\mathcal{A}^{\otimes}/\mathcal{P}^{\otimes}}(X)$
spanned by those objects
$\mathcal{A}^{\otimes}\to X_y$
which are maps of $\infty$-operads over $\mathcal{P}^{\otimes}$.
By restriction,
we obtain a map of $\infty$-categories
\[ \xi_Q: {\rm Alg}_{\mathcal{A}/\mathcal{P}}(X)^{\boxtimes, \vee}
          \longrightarrow \mathcal{Q}^{\boxtimes,\rm op}.\]
\end{definition}

\begin{proposition}\label{prop:A-algebra-Q-monoid}
The map 
$\xi_Q:
{\rm Alg}_{\mathcal{A}/\mathcal{P}}(X)^{\boxtimes, \vee}
\to\mathcal{Q}^{\boxtimes,\rm op}$
is a cartesian fibration of opposite $\infty$-operads.
\end{proposition}

\proof
We set $A={\rm Alg}_{\mathcal{A}/\mathcal{P}}
(X)^{\boxtimes, \vee}$
for simplicity.
It suffices to show that
(1) $\xi_Q$ is a cartesian fibration,
and that
(2) for any $y\in\mathcal{Q}^{\boxtimes}$,
a family of inert morphisms
$\{y\to y_j|\ j\in |J|\}$
induces an equivalence 
$A_y\to \prod_{j\in |J|}A_{y_j}$
of $\infty$-categories.

For (1),
let $a\in A_y$ and 
let $f: y\to y'$ be a morphism of $\mathcal{Q}^{\boxtimes}$.
We take a cartesian lift
$f^*a\to a$ of $f$ in
$F_{\mathcal{A}^{\otimes}/\mathcal{P}^{\otimes}}(X)$.
Since $A$ is a full subcategory
of $F_{\mathcal{A}^{\otimes}/\mathcal{P}^{\otimes}}(X)$,
it suffices to show that $f^*a\in A_{y'}$.
Since $X\to\mathcal{Q}^{\boxtimes,\rm op}$
is a cartesian fibration,
$f$ induces a map
$f^*: X_{y}\to X_{y'}$, which is
a map of $\infty$-operads over $\mathcal{P}^{\otimes}$.
Since $f^*a$ is identified with
the composite $\mathcal{A}^{\otimes}\stackrel{a}{\to}
X_y\stackrel{f^*}{\to}
X_{y'}$,
we see that $f^*a\in A_{y'}$.

For (2),
we note that
there is an equivalence
$A_{y''}\simeq {\rm Alg}_{\mathcal{A}/\mathcal{P}}(X_{y''})$
of $\infty$-categories
for any $y''\in\mathcal{Q}^{\boxtimes}$.
Since the family
$\{y\to y_j|\ j\in |J|\}$ induces an equivalence
$X_y\to
\prod_{j\in |J|}X_{y_j}$
of $\mathcal{P}$-monoidal $\infty$-categories,
we obtain an equivalence
$A_y\to\prod_{j\in |J|}A_{y_j}$.
\qed

%@\bigskip

In a similar manner
we can construct a cocartesian fibration
of $\infty$-operads
${\rm Coalg}_{\mathcal{B}/\mathcal{Q}}(X)^{\otimes}\to
\mathcal{P}^{\otimes}$,
%for an $\infty$-operad $\mathcal{B}^{\boxtimes}$
%over $\mathcal{Q}^{\boxtimes}$,
where objects of
${\rm Coalg}_{\mathcal{B}/\mathcal{Q}}(X)^{\otimes}$
are $\mathcal{B}$-coalgebra objects
in $X_x$
%over $\mathcal{Q}^{\otimes}$
for some $x\in\mathcal{P}^{\otimes}$.
%By Lemma~\ref{lemma:xi-cocartesian-fibration},
%we have a cocartesian fibration
%$F_{\mathcal{B}^{\boxtimes,\rm op}/\mathcal{Q}^{\boxtimes, \rm op}}
%(X)\to\mathcal{P}^{\otimes}$.

\begin{definition}\rm
We define an $\infty$-category
\[ {\rm Coalg}_{\mathcal{B}/\mathcal{Q}}(X)^{\otimes} \]
to be the full subcategory of
$F_{\mathcal{B}^{\boxtimes, \rm op}/\mathcal{Q}^{\boxtimes, \rm op}}
(X)$
spanned by those objects
$\mathcal{B}^{\boxtimes, \rm op}
\to X_x$
which are maps of opposite $\infty$-operads
over $\mathcal{Q}^{\boxtimes, \rm op}$
for some $x\in\mathcal{P}^{\otimes}$.
By restriction,
we obtain a map of $\infty$-categories
\[ \xi_P: {\rm Coalg}_{\mathcal{B}/\mathcal{Q}}(X)^{\otimes}
          \longrightarrow \mathcal{P}^{\otimes}.\]
\end{definition}

We can prove the following proposition
in a similar manner to
Proposition~\ref{prop:A-algebra-Q-monoid}.

\begin{proposition}\label{prop:B-coalgebra-Q-comonoid}
The map 
$\xi_P: {\rm Coalg}_{\mathcal{B}/\mathcal{Q}}(X)^{\otimes}
\to \mathcal{P}^{\otimes}$
is a cocartesian fibration of
$\infty$-operads.
%over $\Phi$.
\end{proposition}

Next,
we will introduce a notion of
bialgebra in a mixed
duoidal $\infty$-category.

\begin{definition}\rm
Let $B: \mathcal{A}^{\otimes}\times\mathcal{B}^{\boxtimes,\rm op}
\to X$
be a map over $\mathcal{P}^{\otimes}\times\mathcal{Q}^{\boxtimes,\rm op}$.
We say that $B$ is an $(\mathcal{A},\mathcal{B})$-bialgebra
of $X$ if the restriction
$B_a: \mathcal{B}^{\boxtimes,\rm op}\to
X_{\overline{a}}$
is a $\mathcal{B}$-coalgebra
%over $Q$
for each $a\in\mathcal{A}^{\otimes}$,
where $\overline{a}\in\mathcal{P}^{\otimes}$
is the image of $a$ under the map
$\mathcal{A}^{\otimes}\to\mathcal{P}^{\otimes}$,
and if the restriction $B_b: \mathcal{A}^{\otimes}\to
X_{\overline{b}}$
is an $\mathcal{A}$-algebra
%over $\mathcal{P}$
for each $b\in \mathcal{B}^{\boxtimes,\rm op}$,
where $\overline{b}\in\mathcal{Q}^{\boxtimes,\rm op}$
is the image of $b$ under the map
$\mathcal{B}^{\boxtimes,\rm op}\to
\mathcal{Q}^{\boxtimes, \rm op}$.
\end{definition}

\begin{definition}\rm
We define an $\infty$-category
\[ \bialgebra_{(\mathcal{A},\mathcal{B})/(\mathcal{P},\mathcal{Q})}(X) \]
to be the full subcategory
of ${\rm Fun}_{\mathcal{P}^{\otimes}\times\mathcal{Q}^{\boxtimes,\rm op}}
(\mathcal{A}^{\otimes}\times\mathcal{B}^{\boxtimes,\rm op},
X)$
spanned by $(\mathcal{A},\mathcal{B})$-bialgebras
of $X$.
%over $(\mathcal{P},\mathcal{Q})$.
When $\mathcal{A}^{\otimes}=\mathcal{P}^{\otimes}$
and $\mathcal{B}^{\boxtimes}=\mathcal{Q}^{\boxtimes}$,
we simply write
$\bialgebra_{(\mathcal{P},\mathcal{Q})}(X)$.
We call ${\rm Bialg}_{(\mathcal{A},\mathcal{B})/(\mathcal{P},\mathcal{Q})}
(X)$ the $\infty$-category of
$(\mathcal{A},\mathcal{B})$-bialgebras
of $X$.
%over $(\mathcal{P},\mathcal{Q})$.
\end{definition}

We will show that
the $\infty$-category of bialgebras
can be described by the $\infty$-categories
of algebras and coalgebras.

\begin{theorem}[{cf.~\cite[Theorem~1.1]{Torii1}}]
\label{thm:bialgebra-algebra-coalgebra}
There are equivalences of $\infty$-categories
\[ {\rm Alg}_{\mathcal{A}/{\mathcal{P}}}
   ({\rm Coalg}_{\mathcal{B}/\mathcal{Q}}
   (X))\simeq
   {\rm BiAlg}_{(\mathcal{A},\mathcal{B})/(\mathcal{P},\mathcal{Q})}
   (X) \simeq
   {\rm Coalg}_{\mathcal{B}/\mathcal{Q}}
   ({\rm Alg}_{\mathcal{A}/{\mathcal{P}}}
   (X)).\]
\end{theorem}

\proof
We will prove the first equivalence.
The second equivalence can be proved in a similar manner.
The $\infty$-category
${\rm Alg}_{\mathcal{A}/\mathcal{P}}
({\rm Coalg}_{\mathcal{B}/\mathcal{Q}}(X))$
is a full subcategory
of
\[ {\rm Fun}_{\mathcal{P}^{\otimes}}
(\mathcal{A}^{\otimes},
F_{\mathcal{B}^{\boxtimes, \rm op}/\mathcal{Q^{\boxtimes, \rm op}}}(X)), \]
which is equivalent to
${\rm Fun}_{\mathcal{P}^{\otimes}\times\mathcal{Q}^{\boxtimes,\rm op}}
(\mathcal{A}^{\otimes}\times\mathcal{B}^{\boxtimes,\rm op},
X)$.
We can verify that
it is spanned by $(\mathcal{A},\mathcal{B})$-bialgebras
of $X$.
%over $(\mathcal{P},\mathcal{Q})$.
\qed

%\newpage
%\input{comodule}
\section{Multiplicative structures
on bicomodules}
\label{section:comodules}

In this section
we study multiplicative structures
on bicomodules over bialgebras.
%Let $\mathcal{O}^{\otimes}$
%be an $\infty$-operad over
%a perfect operator category.
First,
we will construct a generalized
bimodule object in
the $\infty$-category of $\mathcal{O}$-monoidal
$\infty$-categories and lax $\mathcal{O}$-monoidal functors
for an $(\mathcal{O}\times\mathbf{Ass})$-monoidal $\infty$-category
$\mathcal{M}$.
Its underlying object is
the $\infty$-category
${}_A\mathrm{BMod}_B(\mathcal{M})$
of $A$-$B$-bimodules
for $(\mathbf{Ass}\times \mathcal{O})$-algebra
objects $A$ and $B$,
which is left tensored 
over ${}_A\mathrm{BMod}_A(\mathcal{M})$
and right tensored
over ${}_B\mathrm{BMod}_B(\mathcal{M})$.
Next,
we will introduce 
$(A,\Gamma)$-$(B,\Sigma)$-bicomodules
for $(\mathcal{O},\mathbf{Ass})$-bialgebras
$\Gamma$ in 
${}_A\mathrm{BMod}_A(\mathcal{M})$ and
$\Sigma$ in 
${}_B\mathrm{BMod}_B(\mathcal{M})$.
The goal of this section
is to show that
the $\infty$-category
of $(A,\Gamma)$-$(B,\Sigma)$-bicomodules
has the structure of an $\mathcal{O}$-monoidal
$\infty$-category
(Theorem~\ref{thm:monoidal-bicomodules}).

%\newpage

\subsection{The nonsymmetric generalized
  $\infty$-operad $\mathrm{BM}$}

In this subsection
we recall the nonsymmetric generalized
$\infty$-operad $\mathrm{BM}$
whose operadic localization 
is the nonsymmetric $\infty$-operad of bimodules
$\mathbf{BM}^{\otimes}$.

We recall that ${\mathbf\Delta}$ is the simplicial indexing category.
We identify a morphism $\sigma: [m]\to [n]$ in
${\mathbf\Delta}$
with a sequence of integers
$(i_0,\ldots,i_m)$ where $0\le i_0\le \cdots\le i_m\le n$.
For each $k\ge 0$,
the projection $\subdelta{k}\to {\mathbf\Delta}^{\rm op}$
is a cocartesian fibration that is classified
by the functor ${\mathbf\Delta}^{\rm op}\to {\rm Cat}$
given by $[n]\mapsto {\mathbf\Delta}^{\rm op}([k],[n])$.
We can verify that it is a double category.
Thus,
we can regard $\subdelta{k}\to{\mathbf\Delta}^{\rm op}$
as a nonsymmetric generalized $\infty$-operad.

\begin{definition}
[{\rm cf.~\cite[Definitions~7.1.1 and 7.1.6,
and Proposition~7.1.7]{Gepner-Haugseng}}]\rm
We set $\mathrm{BM}=\subdelta{1}$.
We regard the projection $\mathrm{BM}\to{\mathbf\Delta}^{\rm op}$
as a nonsymmetric generalized $\infty$-operad.
Let $\mathbf{BM}^{\otimes}\to{\mathbf\Delta}^{\rm op}$
be the nonsymmetric $\infty$-operad of bimodules.
There is a map $L_{\rm BM}:
\mathrm{BM}\to \mathbf{BM}^{\otimes}$
of nonsymmetric generalized $\infty$-operads
that exhibits $\mathbf{BM}^{\otimes}$ as
the operadic localization of $\mathrm{BM}$.
\end{definition}

There are two full subcategories
of $\mathrm{BM}$ which are isomorphic to
$\mathbf{Ass}^{\otimes}={\mathbf\Delta}^{\rm op}$.
We write $\mathbf{Ass}^{\otimes}_L$ and $\mathbf{Ass}^{\otimes}_R$
for the full subcategories of ${\rm BM}$ spanned by $(0^k)$
and $(1^k)$ for $k>0$,
respectively.

We will introduce algebraic patterns
$\mathrm{BM}^{\flat}$ and $\mathrm{BM}^{\natural}$.
We note that since the map $\mathrm{BM}\to{\mathbf\Delta}^{\rm op}$
is a cocartesian fibration,
we can lift the (inert, active) factorization system
on ${\mathbf\Delta}^{\rm op}$ to $\mathrm{BM}$
by \cite[Proposition~2.1.2.5]{Lurie2}.

\begin{definition}\rm
We define $\mathrm{BM}^{\flat}$ and
$\mathrm{BM}^{\natural}$  
to be algebraic patterns
in which 
$\mathrm{BM}^{\flat,\rm el}$
is spanned by $(0,0),(0,1),(1,1)$,
and
$\mathrm{BM}^{\natural,\rm el}$
is spanned by $(0),(1),(0,0),(0,1),(1,1)$,
respectively.
\end{definition}

Now,
we can consider $\mathrm{BM}$-monoid objects
and generalized $\mathrm{BM}$-monoid objects.

\begin{definition}\rm
Let $\mathcal{C}$ be an $\infty$-category with finite products,
and let $M: \mathrm{BM}\to\mathcal{C}$
be a functor.
We say that $M$
is a $\mathrm{BM}$-monoid object of $\mathcal{C}$
if it is a $\mathrm{BM}^{\flat}$-Segal object.
We denote by
$\mathrm{Mon}_{\mathrm{BM}}(\mathcal{C})$
the full subcategory of
$\mathrm{Fun}(\mathrm{BM},\mathcal{C})$
spanned by $\mathrm{BM}$-monoid objects.

When $\mathcal{C}$ has finite limits,
we say that $M$
is a generalized $\mathrm{BM}$-monoid object
of $\mathcal{C}$
if it is a $\mathrm{BM}^{\natural}$-Segal object.
We denote by $\mathrm{Mon}_{\mathrm{BM}}^{\rm gen}(\mathcal{C})$
the full subcategory of $\mathrm{Fun}(\mathrm{BM},\mathcal{C})$
spanned by generalized $\mathrm{BM}$-objects.
\end{definition}

\begin{remark}\label{remark:two-BM-monoid-objects}
\rm
Since $\mathbf{BM}^{\otimes}$ is the operadic
localization of $\mathrm{BM}$,
the localization functor $L_{BM}$
induces an equivalence of $\infty$-categories
$\mathrm{Mon}_{\mathbf{BM}}(\mathcal{C})\to
\mathrm{Mon}_{\mathrm{BM}}(\mathcal{C})$
for any $\infty$-category $\mathcal{C}$
with finite products. 
\end{remark}

%\newpage

\subsection{Looping construction of bimodules}

In this subsection
we will construct $\mathrm{BM}$-monoid
objects in the $\infty$-category of
$\mathcal{O}$-monoidal $\infty$-categories
and lax $\mathcal{O}$-monoidal functors
from generalized $\mathrm{BM}$-monoid objects
in the $\infty$-category of
$\mathcal{O}$-monoidal $\infty$-categories
and strong $\mathcal{O}$-monoidal functors
by looping construction.

First,
we consider the case in which
$\mathcal{O}^{\otimes}$ is the trivial
symmetric $\infty$-operad.
%$\mathbf{Triv}^{\otimes}$.
Let $\mathcal{C}$ be an $\infty$-category
with finite limits,
and let $M: \mathrm{BM}\to \mathcal{C}$
be a generalized $\mathrm{BM}$-monoid object of $\mathcal{C}$.
Let $A(i): *\to M(i)$ be
morphisms in $\mathcal{C}$ for $i=0,1$,
where $*$ is a final object of $\mathcal{C}$.
We set
\[ N(i,j)=\{A(i)\}\times_{M(i)}M(i,j)
   \times_{M(j)}\{A(j)\}\]
for $i,j=0,1$ with $i\le j$.
We will show that
$N(i,j)$ can be promoted to
a ${\rm BM}$-monoid object. 

\begin{proposition}\label{prop:construction-bimodule}
Let $\mathcal{C}$ be an $\infty$-category with finite limits,
and let $M: \mathrm{BM}\to\mathcal{C}$
be a generalized $\mathrm{BM}$-monoid object of $\mathcal{C}$.
For any morphisms $A(i): *\to M(i)$ for $i=0,1$,
there is a $\mathrm{BM}$-monoid object
${}_{A(0)}M_{A(1)}$
%\in \mathrm{Mon}_{\mathrm{BM}}(\mathcal{C})$
such that
${}_{A(0)}M_{A(1)}(i,j)\simeq N(i,j)$
for $i,j=0,1$ with $i\le j$.
\end{proposition}

\proof
%[Proof of Proposition~\ref{prop:construction-bimodule}]
We let $\mathrm{BM}_{[0]}=\{(0),(1)\}$
be the fiber of the projection at $[0]$.
We denote by $j: \mathrm{BM}_{[0]}\to \mathrm{BM}$
the inclusion.
There is an adjunction
\[ j^*: {\rm Fun}(\mathrm{BM},\mathcal{C})\rightleftarrows
   {\rm Fun}(\mathrm{BM}_{[0]},\mathcal{C}): j_*,\]
where the left adjoint $j^*$ is the restriction,
and the right adjoint $j_*$ is 
the right Kan extension along $j$.
We notice that ${\rm Fun}(\mathrm{BM}_{[0]},\mathcal{C})$
is equivalent to $\mathcal{C}\times\mathcal{C}$ and
that $j^*Z\simeq (Z(0),Z(1))$
for $Z\in {\rm Fun}(\mathrm{BM},\mathcal{C})$.

The morphisms $A(0), A(1)$ determine a map
$(A(0),A(1)): (*,*)\to (M(0),M(1))
   \simeq j^*M$
in $\mathcal{C}\times\mathcal{C}$.
We define ${}_{A(0)}M_{A(1)}$
by the following pullback diagram
\begin{align}
\label{eq:pullback-general-BM-monoid}
  \vcenter{\xymatrix{
     {}_{A(0)}M_{A(1)}\ar[r]\ar[d] & M\ar[d] \\
     j_*(*,*)\ar[r] &
     j_*j^*M \\
   }}\end{align}
in ${\rm Fun}(\mathrm{BM},\mathcal{C})$,
where the bottom horizontal arrow is $j_*(A(0),A(1))$.
We can verify that ${}_{A(0)}M_{A(1)}$
has the desired properties.
\qed

%\newpage
%@\bigskip

We apply Proposition~\ref{prop:construction-bimodule}
to the following situation:
Let $M$
%\in \mathrm{Mon}_{\mathrm{BM}}^{\rm gen}(\mathrm{Mon}_{\mathcal{O}})$
be a generalized $\mathrm{BM}$-monoid
object in $\mathrm{Mon}_{\mathcal{O}}$.
We have the $\infty$-category
$\mathrm{Alg}_{\mathcal{O}}(M(i))$
%$\mathrm{Alg}_{\mathcal{O}}^{\rm Mon}(M(i))$
of $\mathcal{O}$-algebra objects in $M(i)$
for $i=0,1$.
For $A(i)\in \mathrm{Alg}_{\mathcal{O}}(M(i))$,
%For $A(i)\in \mathrm{Alg}_{\mathcal{O}}^{\rm Mon}(M(i))$,
we regard $A(i)$ as a map
$A(i): *\to M(i)$
in ${\rm Op}_{\mathcal{O}}$.
By Proposition~\ref{prop:construction-bimodule},
we obtain a $\mathrm{BM}$-monoid object
${}_{A(0)}M_{A(1)}$ in ${\rm Op}_{\mathcal{O}}$
such that
${}_{A(0)}M_{A(1)}(i,j)\simeq N(i,j)$
for $i,j=0,1$ with $i\le j$.

We would like to have a $\mathrm{BM}$-monoid
object ${}_{A(0)}M_{A(1)}$ in $\mathrm{Mon}_{\mathcal{O}}^{\rm lax}$
instead of ${\rm Op}_{\mathcal{O}}$.
For this purpose,
we have the following theorem.

\begin{theorem}\label{thm:construction-bimodule-O-monoidal}
Let $M$
%\in \mathrm{Mon}_{\mathrm{BM}}^{\rm gen}(\mathrm{Mon}_{\mathcal{O}})$
be a generalized $\mathrm{BM}$-monoid
object of $\mathrm{Mon}_{\mathcal{O}}$,
and let
$A(i)\in \mathrm{Alg}_{\mathcal{O}}(M(i))$
%$A_{(i)}\in \mathrm{Alg}_{\mathcal{O}}^{\rm Mon}(M(i))$
for $i=0,1$.
If 
$M(i,j)^{\otimes}\to M(i)^{\otimes}
\times_{\mathcal{O}^{\otimes}}M(j)^{\otimes}$
is a cocartesian fibration of $\infty$-operads
for $i,j=0,1$ with $i\le j$,
then ${}_{A(0)}M_{A(1)}$
is a $\mathrm{BM}$-monoid object
in $\mathrm{Mon}_{\mathcal{O}}^{\rm lax}$.
\end{theorem}

\proof
We consider pullback diagram~(\ref{eq:pullback-general-BM-monoid})
in the proof of Proposition~\ref{prop:construction-bimodule}.
Since $j_*$ is a right adjoint,
$j_*(*,*)$ is a final object
of ${\rm Fun}(\mathrm{BM},{\rm Op}_{\mathcal{O}})$.
Hence,
in order to show that
${}_{A(0)}M_{A(1)}(0^k,1^l)$
is an $\mathcal{O}$-monoidal $\infty$-category,
it suffices to show that
$M(0^k,1^l)\to j_*j^*M(0^k,1^l)$
is a cocartesian fibration of $\infty$-operads.
Since this map is identified with
a product of copies of the maps
$M(i,j)^{\otimes}\to M(i)^{\otimes}
\times_{\mathcal{O}^{\otimes}}M(j)^{\otimes}$,
the theorem follows from the assumption.
\qed

%\bigskip

%\newpage

\subsection{Multiplicative structures on bi(co)modules}
\label{subsection:multiplicative-bimodules}

In this subsection
we study multiplicative structures
on bi(co)modules.
We will construct an $\infty$-category
of bimodules
%where the ambient $\infty$-category
%is a bimodule object
in a bimodule object of $\monoplax$,
and show that it has the structure of an
$\mathcal{O}$-monoidal $\infty$-category.
By taking opposites,
we also show that an $\infty$-category
of bicomodules has the structure
of an $\mathcal{O}$-monoidal
$\infty$-category.

We recall that 
$\monoplax$
is the $\infty$-category 
of small $\mathcal{O}$-monoidal $\infty$-categories
and oplax $\mathcal{O}$-monoidal functors.
We suppose that
$X\in \mathrm{Mon}_{\mathrm{BM}}(\monoplax)$
is a $\mathrm{BM}$-monoid object of $\monoplax$.
We set
%$X_L=X(0,0),
%X_R=X(1,1)$,
%and
$X_M=X(0,1)$.
By Remarks~\ref{remark:equivalence-mixed-duoidal-bilax}
and \ref{remark:two-BM-monoid-objects},
we can regard $X$
as a mixed $(\mathbf{BM},\mathcal{O})$-duoidal
$\infty$-category
\[ X \longrightarrow \mathbf{BM}^{\otimes}
   \times \mathcal{O}^{\boxtimes,{\rm op}}.\]
By restriction,
we obtain a cartesian fibration of opposite $\infty$-operads
$X_M\to\mathcal{O}^{\boxtimes,{\rm op}}$,
%$X_M^{\boxtimes,\vee}\to
%  \mathcal{O}^{\boxtimes,{\rm op}}$,
and
mixed $(\mathbf{Ass},\mathcal{O})$-duoidal
$\infty$-categories
$X_L\to\mathbf{Ass}^{\otimes}_L\times
   \mathcal{O}^{\boxtimes,{\rm op}}$
and
$X_R\to\mathbf{Ass}^{\otimes}_R\times
   \mathcal{O}^{\boxtimes,{\rm op}}$.

Applying Proposition~\ref{prop:A-algebra-Q-monoid}
to the $(\mathbf{BM},\mathcal{O})$-duoidal $\infty$-category $X$,
%to the map 
%$X\to
%\mathbf{BM}^{\otimes}\times\mathcal{O}^{\boxtimes,\rm op}$,
we obtain a cartesian fibration of opposite $\infty$-operads
\[ p: \mathrm{BMod}(X_M)^{\boxtimes,\vee}\longrightarrow
   \mathcal{O}^{\boxtimes,\rm op}, \]
where $\mathrm{BMod}(X_M)^{\boxtimes,\vee}=
\mathrm{Alg}_{\mathrm{BM}/\mathbf{BM}^{\otimes}}(X)^{\boxtimes,\vee}$.
The map $p$ gives $\mathrm{BMod}(X_M)$
the structure of an $\mathcal{O}$-monoidal $\infty$-category.
Applying Proposition~\ref{prop:A-algebra-Q-monoid}
to the mixed $(\mathbf{Ass},\mathcal{O})$-duoidal
$\infty$-categories $X_L$ and $X_R$,
%$A_v^{\otimes}\to
%\mathbf{Ass}_v^{\otimes}\times\mathcal{O}^{\boxtimes,\rm op}$,
we also obtain a cartesian fibration of opposite $\infty$-operads
\[ q: \mathrm{Alg}(X_L)^{\boxtimes,\vee}
      \times_{\mathcal{O}^{\boxtimes,\rm op}}
      \mathrm{Alg}(X_R)^{\boxtimes,\vee}
      \longrightarrow
   \mathcal{O}^{\boxtimes,\rm op}, \]
where $\mathrm{Alg}(X_v)^{\boxtimes,\vee}=
\mathrm{Alg}_{\mathbf{Ass}_v^{\otimes}/\mathbf{BM}^{\otimes}}
(X)^{\boxtimes,\vee}$
for $v=L,R$.
%The map $q$ gives $\mathrm{Alg}(A_v)$
%the structure of an $\mathcal{O}$-monoidal $\infty$-category.
The naturality of constructions
induces a commutative diagram
\[ \xymatrix{
  \mathrm{BMod}(X_M)^{\boxtimes,\vee}
  \ar[rr]^-{\pi}\ar[dr]_p&&
     \mathrm{Alg}(X_L)^{\boxtimes,\vee}
     \times_{\mathcal{O}^{\boxtimes, \rm op}}
     \mathrm{Alg}(X_R)^{\boxtimes,\vee}
           \ar[dl]^q\\
   &\mathcal{O}^{\boxtimes,{\rm op}}.&\\
  }\]

\begin{proposition}\label{prop:bmod-alg-cartesian}
The map
$\pi: \mathrm{BMod}(X_M)^{\boxtimes,\vee}\to
\mathrm{Alg}(X_L)^{\boxtimes,\vee}
\times_{\mathcal{O}^{\boxtimes,\rm op}}
\mathrm{Alg}(X_R)^{\boxtimes,\vee}$
is a cartesian fibration of opposite $\infty$-operads.
\end{proposition}

%\if0
\proof
By \cite[Proposition~5.45]{Haugseng1},
it suffices to show that
(1) 
$\pi_x:
\mathrm{BMod}(X_M)^{\boxtimes,\vee}_x\to
\mathrm{Alg}(X_L)^{\boxtimes,\vee}_x
\times \mathrm{Alg}(X_R)^{\boxtimes,\vee}_x$
is a cartesian fibration
for any $x\in \mathcal{O}^{\boxtimes,{\rm op}}$,
(2)
$\pi:
\mathrm{BMod}(X_M)^{\boxtimes,\vee}\to
\mathrm{Alg}(X_L)^{\boxtimes,\vee}
\times_{\mathcal{O}^{\boxtimes,\rm op}}
\mathrm{Alg}(X_R)^{\boxtimes,\vee}$
carries $p$-cartesian edges
to $q$-cartesian edges, and
(3)
for any map $\phi: x\to y$ in $\mathcal{O}^{\boxtimes,{\rm op}}$,
the induced map
$\phi^*: \mathrm{BMod}(X_M)^{\boxtimes,\vee}_y
\to \mathrm{BMod}(X_M)^{\boxtimes,\vee}_x$
carries $\pi_y$-cartesian edges
to $\pi_x$-cartesian edges.

For (1),
the map $\pi_x$
is equivalent to the projection map
$\mathrm{BMod}(M_x)\to \mathrm{Alg}(X_{L,x})
\times\mathrm{Alg}(X_{R,x})$
which is a cartesian fibration.
For (2),
we let $f: (a,m,b)\to (c,n,d)$ be a cartesian edge of
$\mathrm{BMod}(X_M)^{\boxtimes,\vee}$
over $\phi: x\to y$ in $\mathcal{O}^{\boxtimes,{\rm op}}$.
We may assume that $x\in \mathcal{O}$
and $\phi: y\to x$ is active in $\mathcal{O}^{\boxtimes}$.
We suppose that $y\in \mathcal{O}^{\boxtimes}_I$
and take a set of inert morphisms
$\{y\to y_i|\ i\in |I|\}$
over $\{\rho_i: I\to \{i\}|\ i\in |I|\}$.
Since $f$ is a cartesian edge,
it induces an equivalence
$(a,m,b)\stackrel{\simeq}{\to}
\boxtimes_{\phi}(c_i,n_i,d_i)$
in $\mathrm{BMod}(X_M)^{\boxtimes,\vee}_x$,
where $(c_i,n_i,d_i)\to (c,n,d)$
is a cartesian edge over $y_i\to y$.
This implies that there is an equivalence
$a\stackrel{\simeq}{\to}\boxtimes_{\phi}(c_i)$
in $\mathrm{Alg}(X_L)^{\boxtimes,\vee}_x$,
and
$b\stackrel{\simeq}{\to}\boxtimes_{\phi}(d_i)$
in $\mathrm{Alg}(X_R)^{\boxtimes,\vee}_x$.
This means that $\pi$
preserves cartesian edges.
For (3),
we may assume that $x\in\mathcal{O}$ and
$\phi: y\to x$ is active in $\mathcal{O}^{\boxtimes}$.
Then,
we have equivalences
$\mathrm{BMod}(X_M)^{\boxtimes,\vee}_y\simeq
\prod_{i\in |I|}\mathrm{BMod}((X_M)_{y_i})$
and
$\mathrm{BMod}(X_M)^{\boxtimes,\vee}_x\simeq
\mathrm{BMod}((X_M)_x)$.
Let $(f,h,g): (a,m,b)\to (c,n,d)$
be a $\pi_y$-cartesian edge of
$\mathrm{BMod}(X_M)^{\boxtimes,\vee}_y$.
Under the equivalence 
$\mathrm{BMod}(X_M)^{\boxtimes,\vee}_y\simeq
\prod_{i\in |I|}\mathrm{BMod}((X_M)_{{y_i}})$,
we regard $(f,h, g)$ as a product of $(f_i,h_i,g_i)$,
where $(f_i,h_i,g_i): (a_i,m_i,b_i)\to (c_i,n_i,d_i)$
is a $\pi_{y_i}$-cartesian edge
in $\mathrm{BMod}(X_M)^{\boxtimes,\vee}_{y_i}
\simeq \mathrm{BMod}((X_M)_{y_i})$.
This implies that $(f_i,h_i,g_i)$
induces an equivalence
$m_i\stackrel{\simeq}{\to} f_i^*n_ig_i^*$
in ${}_{a_i}\mathrm{BMod}_{b_i}((X_M)_{y_i})$.
The functor $\phi^*$ takes $(f,h,g)$
to $\boxtimes_{\phi}(f_i,h_i,g_i):
\boxtimes_{\phi}(a_i,m_i,b_i)\to
\boxtimes_{\phi}(c_i,n_i,d_i)$.
Since $\boxtimes_{\phi}(f_i,h_i,g_i)$
induces an equivalence
$\boxtimes_{\phi}h_i: \boxtimes_{\phi}(m_i)
\stackrel{\simeq}{\to}(\boxtimes_{\phi}(f_i))^*
(\boxtimes_{\phi}(n_i))(\boxtimes_{\phi}(g_i))^*$,
we see that 
$\boxtimes_{\phi}(f_i,h_i,g_i)$ is
a $\pi_x$-cartesian edge.
\qed
%\fi

%\bigskip

Let $\Gamma(v)\in \bialgebra_{(\mathbf{Ass},\mathcal{O})}(X_v)$
be an $(\mathbf{Ass},\mathcal{O})$-bialgebra object of
the mixed $(\mathbf{Ass},\mathcal{O})$-duoidal
$\infty$-category $X_v$ for $v=L,R$.
By Theorem~\ref{thm:bialgebra-algebra-coalgebra},
this determines a map
\[ (\Gamma(L),\Gamma(R)):
   \mathcal{O}^{\boxtimes,{\rm op}}\longrightarrow
   \mathrm{Alg}(X_L)^{\boxtimes,\vee}
   \times_{\mathcal{O}^{\boxtimes,\rm op}}
   \mathrm{Alg}(X_R)^{\boxtimes,\vee}\]
of opposite $\infty$-operads over $\mathcal{O}^{\boxtimes,\rm op}$.
We consider the following pullback diagram
\begin{equation}\label{eq:diagram-bmod-pullback}
\vcenter{  \xymatrix{
    {}_{\Gamma(L)}\mathrm{BMod}_{\Gamma(R)}(X_M)^{\boxtimes,\vee}
    \ar[r]\ar[d]&
    \mathrm{BMod}(X_M)^{\boxtimes,\vee}\ar[d]^{\pi}\\
    \mathcal{O}^{\boxtimes,{\rm op}}\ar[r]^-{(\Gamma(L),\Gamma(R))}&
    \mathrm{Alg}(X_L)^{\boxtimes,\vee}
    \times_{\mathcal{O}^{\boxtimes,\rm op}}
    \mathrm{Alg}(X_R)^{\boxtimes,\vee}.\\
   }}\end{equation}
By Proposition~\ref{prop:bmod-alg-cartesian},
we obtain the following theorem.

\begin{theorem}\label{thm:Gamma-Bmod-O-monoidal}
Let $X$ be a $\mathrm{BM}$-monoid object
of $\monoplax$.
For $(\mathbf{Ass},\mathcal{O})$-bialgebra objects
$\Gamma(v)$ of
the mixed $(\mathbf{Ass},\mathcal{O})$-duoidal
$\infty$-categories $X_v$ for $v=L,R$,
the map
${}_{\Gamma(L)}\mathrm{BMod}_{\Gamma(R)}(X_M)^{\boxtimes,\vee}\to
\mathcal{O}^{\boxtimes,{\rm op}}$ is a cartesian fibration
of opposite $\infty$-operads,
and hence it gives
${}_{\Gamma(L)}\mathrm{BMod}_{\Gamma(R)}(X_M)$
the structure of an $\mathcal{O}$-monoidal $\infty$-category.
Furthermore,
the forgetful functor
${}_{\Gamma(L)}\mathrm{BMod}_{\Gamma(R)}(X_M)\to X_M$
is strong $\mathcal{O}$-monoidal.
\end{theorem}

\proof
By Proposition~\ref{prop:bmod-alg-cartesian},
the right vertical arrow $\pi$ in
(\ref{eq:diagram-bmod-pullback})
is a cartesian fibration.
Thus,
so is the left vertical arrow.
%By its proof,
We see that the forgetful functor
preserves cartesian edges
over $\mathcal{O}^{\boxtimes,\rm op}$.
\qed

Next,
we consider a dual situation of the above.
Let $Y$ be a $\mathrm{BM}$-monoid object
of $\monlax$.
We regard it as 
a mixed $(\mathcal{O},\mathbf{BM})$-duoidal
$\infty$-category
$Y\to \mathcal{O}^{\boxtimes}\times\mathbf{BM}^{\otimes,\rm op}$
by Remark~\ref{remark:equivalence-mixed-duoidal-bilax}.
By Proposition~\ref{prop:duality-mixed-duoidal},
the opposite $Y^{\rm op}\to
\mathbf{BM}^{\otimes}\times \mathcal{O}^{\boxtimes,\rm op}$
has the structure of a mixed $(\mathbf{BM},\mathcal{O})$-duoidal
$\infty$-category.
We obtain a cartesian fibration
${}_{\Gamma(L)}\mathrm{BMod}_{\Gamma(R)}(Y_M^{\rm op})^{\boxtimes,\vee}\to
\mathcal{O}^{\boxtimes,\rm op}$
of opposite $\infty$-operads
by Theorem~\ref{thm:Gamma-Bmod-O-monoidal},
where 
$\Gamma(L)$ and $\Gamma(R)$ are
$(\mathbf{Ass},\mathcal{O})$-bialgebra objects
of the mixed $(\mathbf{Ass},\mathcal{O})$-duoidal
$\infty$-categories $Y_L^{\rm op}$ and $Y_R^{\rm op}$,
respectively.
We can regard $\Gamma(L)$ and $\Gamma(R)$ as
$(\mathcal{O},\mathbf{Ass})$-bialgebra objects
of the mixed $(\mathcal{O},\mathbf{Ass})$-duoidal
$\infty$-categories
$Y_L$ and $Y_R$,
respectively.
Roughly speaking,
$\Gamma(v)$ is an object of $Y_v$ for $v=L,R$
equipped with
an associative coproduct and a counit
\[ \Gamma(v)_{x}\longrightarrow
   \Gamma(v)_{x}\otimes_{x}\Gamma(v)_{x},
   \qquad \Gamma(v)_{x}\to 1_x \]
for each $x\in\mathcal{O}^{\otimes}$,
and an $\mathcal{O}$-monoidal product
\[ \boxtimes_{\phi}: \Gamma(v)_{x}
   \longrightarrow \Gamma(v)_{y} \]
for each active morphism $\phi: x\to y$
of $\mathcal{O}^{\otimes}$
that are compatible in suitable sense.

\begin{definition}\rm
We define a cocartesian fibration of
$\infty$-operads
\[ {}_{\Gamma(L)}\bicomod_{\Gamma(R)}(Y_M)^{\boxtimes}
   \to\mathcal{O}^{\boxtimes}, \] 
to be the opposite of
the map
${}_{\Gamma(L)}\mathrm{BMod}_{\Gamma(R)}(Y_M^{\rm op})^{\boxtimes,\vee}\to
\mathcal{O}^{\boxtimes,\rm op}$.
We call 
${}_{\Gamma(L)}\bicomod_{\Gamma(R)}(Y_M)$
the $\infty$-category of 
$\Gamma(L)$-$\Gamma(R)$-bicomodules of $Y_M$.
\end{definition}

Roughly speaking,
an object of
${}_{\Gamma(L)}\bicomod_{\Gamma(R)}(Y_M)$
is an object $M_x$ of $(Y_M)_x$ for some $x\in\mathcal{O}$
equipped with a coaction map
\[ M_x\longrightarrow \Gamma(L)_{x}\otimes_{x}
   M_x\otimes_{x} \Gamma(R)_{x} \]
that is coassociative and counital
up to higher coherent homotopy.   

In summary we obtain the following corollary.

\begin{corollary}\label{cor:mixed-to-product-comodules}
Let $Y$ be a $\mathrm{BM}$-monoid object
of $\monlax$.
For $(\mathcal{O},\mathbf{Ass})$-bialgebra objects
$\Gamma(L)$ and $\Gamma(R)$
of the mixed $(\mathcal{O},\mathbf{Ass})$-duoidal
$\infty$-categories $Y_L$ and $Y_R$,
respectively,
the $\infty$-category ${}_{\Gamma(L)}\bicomod_{\Gamma(R)}(Y_M)$
of $\Gamma(L)$-$\Gamma(R)$-bicomodules in $Y_M$
has the structure of an $\mathcal{O}$-monoidal
$\infty$-category such that
the forgetful functor
${}_{\Gamma(L)}\bicomod_{\Gamma(R)}(Y_M)\to Y_M$
is strong $\mathcal{O}$-monoidal.
\end{corollary}

\subsection{Bicomodules in
$(\mathcal{O}\times\mathbf{Ass})$-monoidal
$\infty$-categories}

In this subsection
we start
with an $(\mathcal{O}\times\mathbf{Ass})$-monoidal
$\infty$-category $\mathcal{M}$ and
$(\mathcal{O}\times\mathbf{Ass})$-algebras $A, B$ in $\mathcal{M}$,
and
construct a $\mathrm{BM}$-monoid object
${}_A\mathbb{BM}_B(\mathcal{M})$ in $\monlax$
such that ${}_A\mathbb{BM}_B(\mathcal{M})(0,1)\simeq
{}_A\mathrm{BMod}_B(\mathcal{M})$.
By
%using Remark~\ref{remark:two-BM-monoid-objects} and 
%applying Theorem~\ref{thm:Gamma-Bmod-O-monoidal}
applying Corollary~\ref{cor:mixed-to-product-comodules}
to ${}_A\mathbb{BM}_B(\mathcal{M})$,
%to the composite of ${}_A\mathbb{BM}_B$ with
%the equivalence $(-)^{\rm op}: \monlax\to\monoplax$,
we obtain an $\mathcal{O}$-monoidal structure
on the $\infty$-category of bicomodules in $\mathcal{M}$.

In order to construct the desired $\mathrm{BM}$-monoid
object ${}_A\mathbb{BM}_B(\mathcal{M})$,
we need a good tensor product
with respect to geometric realizations.
We denote by
$\cat^{\rm geo}$
the subcategory of $\cat$
spanned by $\infty$-categories
which have geometric realizations,
and functors which preserve them.
Since $\cat^{\rm geo}$ has finite products,
we can consider the $\infty$-category
$\mons^{\rm geo}=\mons(\cat^{\rm geo})$ of
$\mathcal{O}$-monoid objects in $\cat^{\rm geo}$.
We can also consider
$\mons^{\rm lax, geo}=\monlax(\cat^{\rm geo})$
the $\infty$-category of $\mathcal{O}$-monoidal
$\infty$-categories 
and lax $\mathcal{O}$-monoidal functors
in $\cat^{\rm geo}$.

We apply Theorem~\ref{thm:construction-bimodule-O-monoidal}
to the following situation:
Let $\mathcal{M}$ be an $(\mathcal{O}\times\mathbf{Ass})$-monoid
object of $\cat^{\rm geo}$.
Roughly speaking,
the $\infty$-category $\mathcal{M}$
has two monoidal tensor products
$\otimes$ and $\boxtimes$,
where $\otimes$ gives
an $\mathbf{Ass}$-monoidal structure
and $\boxtimes$ gives an $\mathcal{O}$-monoidal structure.
They are compatible in the sense that
the functors
\[ \otimes: \mathcal{M}\times\mathcal{M}\to \mathcal{M},
    \qquad 1_{\otimes}: [0]\to \mathcal{M}\]
are strong $\mathcal{O}$-monoidal.
Furthermore,
the $\infty$-category
$\mathcal{M}_x$ has geometric realizations
for each $x\in\mathcal{O}$,
and the monoidal tensor products $\otimes$ and $\boxtimes$
preserve geometric realizations
separately in each variable.

We regard $\mathcal{M}$ as an object
of $\mathrm{Mon}_{\mathbf{Ass}}(\mathrm{Mon}_{\mathcal{O}}^{\rm geo})$
by Lemma~\ref{lemma:monoidal-category-over-tensor-product}.
We can construct a category object
$\mathrm{ALG}(\mathcal{M})\in
\mathrm{Mon}_{\mathbf{Ass}}^{\rm gen}
(\mathrm{Mon}_{\mathcal{O}}^{\rm geo})$
such that
$\mathrm{ALG}(\mathcal{M})([0])\simeq
\mathrm{Alg}(\mathcal{M})$ and
$\mathrm{ALG}(\mathcal{M})([1])\simeq
\mathrm{BMod}(\mathcal{M})$
by \cite[Lemma~4.19 and Theorem~4.39]{Haugseng1}.
By taking a pullback
%of the category object $\mathrm{ALG}(\mathcal{M})$
along the projection map
$\mathrm{BM}\to\mathbf{\Delta}^{\rm op}=\mathbf{Ass}^{\otimes}$,
we obtain a generalized $\mathrm{BM}$-monoid
object $\mathbb{BM}(\mathcal{M})$
of $\mathrm{Mon}_{\mathcal{O}}^{\rm geo}$
such that $\mathbb{BM}(\mathcal{M})(i,j)
\simeq \mathrm{BMod}(\mathcal{M})$
and $\mathbb{BM}(\mathcal{M})(i)\simeq
\mathrm{Alg}(\mathcal{M})$.
We notice that
there is an equivalence 
$\mathrm{Alg}_{\mathcal{O}}
(\mathrm{Alg}(\mathcal{M}))
\simeq
\mathrm{Alg}_{(\mathcal{O}\times\mathbf{Ass})}(\mathcal{M})$
of $\infty$-categories
by Proposition~\ref{prop:double-algebra-decomposition}.

\begin{proposition}\label{prop:construction-BMod-object}
Let $\mathcal{M}$
be an $(\mathcal{O}\times\mathbf{Ass})$-monoid
object of $\cat^{\rm geo}$.
For any $(\mathcal{O}\times\mathbf{Ass})$-algebra
objects $A,B$ in $\mathcal{M}$,
there is a $\mathrm{BM}$-monoid object
${}_A\mathbb{BM}_B(\mathcal{M})$
in $\mathrm{Mon}_{\mathcal{O}}^{\rm lax, geo}$
such that 
\[ {}_A\mathbb{BM}_B(\mathcal{M})(0,0)\simeq
   {}_A\mathrm{BMod}_A(\mathcal{M}),\quad
   {}_A\mathbb{BM}_B(\mathcal{M})(1,1)\simeq
   {}_B\mathrm{BMod}_B(\mathcal{M}), \]
\[ {}_A\mathbb{BM}_B(\mathcal{M})(0,1)\simeq
   {}_A\mathrm{BMod}_B(\mathcal{M}).\]
%${}_A\mathbb{BM}_B(\mathcal{M})(0,0)\simeq{}_A
%   \mathrm{BMod}_A(\mathcal{M})$,
%${}_A\mathbb{BM}_B(\mathcal{M})(0,1)\simeq
%          {}_A\mathrm{BMod}_B(\mathcal{M})$,
%          and
%${}_A\mathbb{BM}_B(\mathcal{M})(1,1)\simeq {}_B
%   \mathrm{BMod}_B(\mathcal{M})$.
\end{proposition}

In order to prove Proposition~\ref{prop:construction-BMod-object},
we make some preliminaries.
We denote by 
\[ p: \mathrm{BMod}(\mathcal{M})^{\boxtimes}\longrightarrow
   \mathcal{O}^{\boxtimes} \]
a cocartesian fibration of $\infty$-operads
which gives
$\mathrm{ALG}(\mathcal{M})([1])\simeq
\mathrm{BMod}(\mathcal{M})$
the structure of the $\mathcal{O}$-monoidal $\infty$-category.
We also denote by 
\[ q: \mathrm{Alg}(\mathcal{M})^{\boxtimes}\longrightarrow
   \mathcal{O}^{\boxtimes} \]
a cocartesian fibration of $\infty$-operads
which gives $\mathrm{ALG}(\mathcal{M})([0])
\simeq\mathrm{Alg}(\mathcal{M})$
the structure of the $\mathcal{O}$-monoidal
$\infty$-category.
We have a commutative diagram
\[ \xymatrix{
  \mathrm{BMod}(\mathcal{M})^{\boxtimes}
  \ar[rr]^-{\pi}\ar[dr]_p&&
     \mathrm{Alg}(\mathcal{M})^{\boxtimes}
     \times_{\mathcal{O}^{\boxtimes}}
     \mathrm{Alg}(\mathcal{M})^{\boxtimes}
           \ar[dl]^{q\times q}\\
   &\mathcal{O}^{\boxtimes},&\\
  }\]
where $\pi$
is induced by the inert morphisms
$[0]\to [1]$ in ${\mathbf\Delta}$.

The following lemma can be proved
by a dual argument to the proof of
Proposition~\ref{prop:bmod-alg-cartesian}
by using the fact that
the tensor product $\boxtimes$
preserves geometric realizations
separately in each variable.

\begin{lemma}\label{lemma:BMod-Alg-cocartesian}
The map
$\pi: \mathrm{BMod}(\mathcal{M})^{\boxtimes}\to
\mathrm{Alg}(\mathcal{M})^{\boxtimes}
\times_{\mathcal{O}^{\boxtimes}}
\mathrm{Alg}(\mathcal{M})^{\boxtimes}$
is a cocartesian fibration of $\infty$-operads.
\end{lemma}

\proof[Proof of Proposition~\ref{prop:construction-BMod-object}]
This follows from
Theorem~\ref{thm:construction-bimodule-O-monoidal}
and Lemma~\ref{lemma:BMod-Alg-cocartesian}.
\qed

%\newpage

%\bigskip

We can regard the $\mathrm{BM}$-monoid
object ${}_A\mathbb{BM}_B(\mathcal{M})$ in $\mons^{\rm lax, geo}$
as a mixed $(\mathcal{O},\mathbf{BM})$-duoidal
$\infty$-category
by Remark~\ref{remark:equivalence-mixed-duoidal-bilax}.
Let $\Gamma$ and $\Sigma$
be $(\mathcal{O},\mathbf{Ass})$-bialgebra objects
of the $(\mathcal{O},\mathbf{Ass})$-duoidal
$\infty$-categories ${}_A\mathbb{BM}_B(\mathcal{M})_L$ and
${}_A\mathbb{BM}_B(\mathcal{M})_R$,
respectively.
We notice that
%since
%${}_A\mathbb{BM}_B(0,0)\simeq
%{}_A\mathrm{BMod}_A(\mathcal{M})$
%and
%${}_A\mathbb{BM}_B(1,1)\simeq
%{}_B\mathrm{BMod}_B(\mathcal{M})$,
$\Gamma$ and $\Sigma$
correspond to $(\mathcal{O},\mathbf{Ass})$-bialgebras
of ${}_A\mathrm{BMod}_A(\mathcal{M})$
and ${}_B\mathrm{BMod}_B(\mathcal{M})$,
respectively.

\begin{definition}\rm
We write 
\[ {}_{(A,\Gamma)}\bicomod_{(B,\Sigma)}
   (\mathcal{M})^{\boxtimes}
   \longrightarrow \mathcal{O}^{\boxtimes} \]
for the cocartesian fibration of
$\infty$-operads
${}_{\Gamma}\mathrm{BMod}_{\Sigma}
({}_A\mathbb{BM}_B(\mathcal{M})_M)^{\boxtimes}
\to \mathcal{O}^{\boxtimes}$,
and call it
%${}_{(A,\Gamma)}\bicomod_{(B,\Sigma)}
%(\mathcal{M})$
the $\infty$-category of 
$(A,\Gamma)$-$(B,\Sigma)$-bicomodules of $\mathcal{M}$.
\end{definition}

By applying Corollary~\ref{cor:mixed-to-product-comodules}
to ${}_A\mathbb{BM}_B(\mathcal{M})$,
we obtain the following theorem.

\begin{theorem}\label{thm:monoidal-bicomodules}
Let $\mathcal{M}$ be
an $(\mathcal{O}\times\mathbf{Ass})$-monoid
object of $\cat^{\rm geo}$, and
let $A, B$ be $(\mathcal{O}\times\mathbf{Ass})$-algebra
objects of $\mathcal{M}$.
For $(\mathcal{O}, \mathbf{Ass})$-bialgebras
$\Gamma$ in ${}_A\mathrm{BMod}_A(\mathcal{M})$
and
$\Sigma$ in ${}_B\mathrm{BMod}_B(\mathcal{M})$,
the $\infty$-category
${}_{(A,\Gamma)}\bicomod_{(A,\Sigma)}(\mathcal{M})$
of $(A,\Gamma)$-$(B,\Sigma)$-bicomodules of $\mathcal{M}$
has the structure of an $\mathcal{O}$-monoidal
$\infty$-category in $\cat^{\rm geo}$ such that 
the forgetful functor
${}_{(A,\Gamma_L)}\bicomod_{(B,\Gamma_R)}(\mathcal{M})\to
{}_A{\rm BMod}_B(\mathcal{M})$
is strong $\mathcal{O}$-monoidal.
%in $\cat^{\rm geo}$.
\end{theorem}

\if0
\begin{remark}\rm
By \cite[Corollary 4.3.3.3]{Lurie2},
the $\mathcal{O}$-monoidal $\infty$-category
${}_{(A,\Gamma)}\bicomod_{(B,\Sigma)}
(\mathcal{M})$
is an object of $\mathrm{Mon}_{\mathcal{O}}^{\rm geo}$,
and the forgetful functor
${}_{(A,\Gamma)}\bicomod_{(B,\Sigma)}
(\mathcal{M})\to {}_A\mathrm{BMod}_B(\mathcal{M})$
is a morphism of $\mathrm{Mon}_{\mathcal{O}}^{\rm geo}$.
\end{remark}
\fi

For an active morphism
$\phi: x\to y$ in $\mathcal{O}^{\boxtimes}$
with $y\in\mathcal{O}$,
we denote by
$({}_A\boxtimes_B)_{\phi}:
{}_A\mathrm{BMod}_B(\mathcal{M})_x\to
{}_A\mathrm{BMod}_B(\mathcal{M})_y$
the tensor product
associated to $\phi$.
For $M_x\in {}_{(A,\Gamma)}\bicomod_{(B,\Sigma)}(\mathcal{M})_x$,
roughly speaking,
the bicomodule structure
on the tensor product
$({}_A\boxtimes_{B})_{\phi} (M_x)$ is
given by
\[ \begin{array}{rcl}
  ({}_A\boxtimes_{B})_{\phi}(M_x)&\longrightarrow&
  ({}_A\boxtimes_{B})_{\phi}
  (\Gamma_x\otimes_{A_x}
  M_x\otimes_{B_x}\Sigma_{x})\\[2mm]
  &\longrightarrow&
  ({}_A\boxtimes_A)_{\phi}(\Gamma_x)
  \otimes_{A_y}
  ({}_A\boxtimes_B)_{\phi}(M_x)
  \otimes_{B_y}
  ({}_B\boxtimes_B)_{\phi}(\Sigma_x)\\[2mm]
  &\longrightarrow&
  \Gamma_y\otimes_{A_y}
  ({}_A\boxtimes_B)_{\phi}(M_x)
  \otimes_{B_y}\Sigma_y,\\
\end{array}\]
where the first arrow is induced by
the bicomodule structure on $M_x$,
the second arrow is an interchange law
of the $\mathrm{BM}$-monoid object
${}_A\mathbb{BM}_B(\mathcal{M})$
in $\mons^{\rm lax, geo}$
in Proposition~\ref{prop:construction-BMod-object},
and
the third arrow is given by
the multiplication maps
of $\Gamma$ and $\Sigma$.

An important case is
when $\mathcal{O}^{\boxtimes}$
is the product
$(\mathbf{Ass}^{\otimes})^n$
of $n$-copies of $\mathbf{Ass}^{\otimes}$
which is an $\infty$-operad
over the perfect operator category
$({\mathbf\Delta}^{\rm op})^n$.
The symmetrization of 
$(\mathbf{Ass}^{\otimes})^n$ 
is equivalent to the little $n$-cubes operad
$\mathbb{E}_n^{\otimes}$,
and hence $\mathbf{Ass}^n$-algebra objects
describe
$\mathbb{E}_n$-algebra objects
in a symmetric $\infty$-operad
(cf.~\cite[Appendix~A]{Haugseng1}).

\begin{corollary}
Let $\mathcal{M}$ be an
$\mathbf{Ass}^{n+1}$-monoid
object of $\cat^{\rm geo}$
and
let $A, B$ be $\mathbf{Ass}^{n+1}$-algebra
objects of $\mathcal{M}$.
For $(\mathbf{Ass}^n, \mathbf{Ass})$-bialgebras
$\Gamma$ of ${}_A\mathrm{BMod}_A(\mathcal{M})$
and $\Sigma$ of ${}_B\mathrm{BMod}_B(\mathcal{M})$,
the $\infty$-category
${}_{(A,\Gamma)}\bicomod_{(B,\Sigma)}(\mathcal{M})$
of $(A,\Gamma)$-$(B,\Sigma)$-bicomodules of $\mathcal{M}$
has the structure of an $\mathbf{Ass}^n$-monoidal
$\infty$-category in $\cat^{\rm geo}$ such that
the forgetful functor
${}_{(A,\Gamma)}\bicomod_{(B,\Sigma)}(\mathcal{M})\to
{}_A\mathrm{BMod}_B(\mathcal{M})$
is strong $\mathbf{Ass}^n$-monoidal.
%in $\cat^{\rm geo}$.
\end{corollary}

%\newpage
%\input{algebra-higher-monoidal}

%\newpage
%\input{looping}
\section{Multiplicative structures
  on right comodules}
\label{section:right-comodules}

In this section
we study multiplicative structures
on right comodules over bialgebras.
The left comodule case can be treated by duality.
Let $\mathcal{M}$ be an
$(\mathbf{Ass}\times\mathcal{O})$-monoid
object of $\cat^{\rm geo}$,
and let $A$ be an $(\mathbf{Ass}\times \mathcal{O})$-monoid
object of $\mathcal{M}$.
First,
we will construct an $\mathrm{RM}$-monoid object
of ${\rm Mon}_{\mathcal{O}}^{\rm lax, geo}$.
Its underlying object is
${\rm RMod}_A(\mathcal{M})$
which is right tensored over
${}_A{\rm BMod}_A(\mathcal{M})$.
Furthermore,
we will introduce right $(A,\Gamma)$-comodules
for $(\mathcal{O},\mathbf{Ass})$-bialgebra
$\Gamma$ in ${}_A{\rm BMod}_A(\mathcal{M})$.
The goal of this section is to show that
the $\infty$-category of right $(A,\Gamma)$-comodules
has the structure of an $\mathcal{O}$-monoidal
$\infty$-category
(Theorem~\ref{thm:right-comodule-main-theorem}).

%\newpage

\subsection{An approximation of
the $\infty$-operad of right modules}

In this subsection
we recall a nonsymmetric generalized $\infty$-operad
$\mathrm{RM}$ 
whose operadic localization 
is the nonsymmetric $\infty$-operad
$\mathbf{RM}^{\otimes}$ of right modules.
After that,
we introduce $\mathrm{RM}$-monoid objects
and generalized $\mathrm{RM}$-monoid objects.
%of an $\infty$-category with finite limits.

\begin{definition}
[{\rm cf.~\cite[Definitions~7.1.1 and 7.1.6,
and Proposition~7.1.7]{Gepner-Haugseng}}]\rm
We let $\mathrm{RM}$ be the full subcategory of
$\mathrm{BM}$ spanned by objects $(0,1^n)$ and $(1,1^n)$
for $n\ge 0$.
The projection $\mathrm{RM}\to\mathbf{\Delta}^{\rm op}$
gives it the structure of
a nonsymmetric generalized $\infty$-operad.
Let $\mathbf{RM}^{\otimes}\to\mathbf{\Delta}^{\rm op}$
be the nonsymmetric $\infty$-operad of right modules.
There is a map of nonsymmetric generalized $\infty$-operads
$L_{\mathrm{RM}}: \mathrm{RM}\to \mathbf{RM}^{\otimes}$
which exhibits $\mathbf{RM}^{\otimes}$ as
the operadic localization of $\mathrm{RM}$.
%The full subcategory $\mathbf{Ass}^{\otimes}_R$
%of $\mathrm{BM}$ is contained in 
%$\mathrm{RM}$.
\end{definition}

We will introduce algebraic patterns
$\mathrm{RM}^{\flat}$ and $\mathrm{RM}^{\natural}$.
We note that since the map $\mathrm{RM}\to{\mathbf\Delta}^{\rm op}$
is a cocartesian fibration,
we can lift the (inert, active) factorization system
on ${\mathbf\Delta}^{\rm op}$ to $\mathrm{RM}$
by \cite[Proposition~2.1.2.5]{Lurie2}.
%By \cite[Definition~9.3, Lemma~9.4,
%Definition~9.9, and Lemma~9.10]{CH},

\begin{definition}\rm
We define $\mathrm{RM}^{\flat}$ and
$\mathrm{RM}^{\natural}$  
to be algebraic patterns
in which $\mathrm{RM}^{\flat,\rm el}$
is spanned by $(0,1),(1,1)$,
and $\mathrm{RM}^{\natural,\rm el}$
is spanned by $(1),(0,1),(1,1)$,
respectively.
\end{definition}

We introduce $\mathrm{RM}$-monoid objects
and generalized $\mathrm{RM}$-monoid objects.

\begin{definition}\rm
Let $\mathcal{C}$ be an $\infty$-category with finite products,
and let $X: \mathrm{RM}\to\mathcal{C}$
be a functor.
We say that $X$
is an $\mathrm{RM}$-monoid object of $\mathcal{C}$
if it is an $\mathrm{RM}^{\flat}$-Segal object.
We denote by $\mathrm{Mon}_{\mathrm{RM}}(\mathcal{C})$
the full subcategory of
$\mathrm{Fun}(\mathrm{RM},\mathcal{C})$
spanned by $\mathrm{RM}$-monoid objects.

When $\mathcal{C}$ has finite limits,
we say that $X$ is 
a generalized $\mathrm{RM}$-monoid object of $\mathcal{C}$
if it is an $\mathrm{RM}^{\natural}$-Segal object.
We denote by $\mathrm{Mon}_{\mathrm{RM}}^{\rm gen}(\mathcal{C})$
the full subcategory of
$\mathrm{Fun}(\mathrm{RM},\mathcal{C})$
spanned by generalized $\mathrm{RM}$-monoid objects.
\end{definition}

\begin{remark}\label{remark:two-RM-monoid-objects}
\rm
As in the case of the $\infty$-operad
of bimodules,
since $\mathbf{RM}^{\otimes}$ is the operadic
localization of $\mathrm{RM}$,
the localization functor $L_{RM}$
induces an equivalence of $\infty$-categories
$\mathrm{Mon}_{\mathbf{RM}}(\mathcal{C})\to
\mathrm{Mon}_{\mathrm{RM}}(\mathcal{C})$
for any $\infty$-category $\mathcal{C}$
with finite products. 
\end{remark}

%\bigskip

Now,
we study a relationship between
(generalized) $\mathrm{RM}$-monoid objects
and (generalized) $\mathrm{BM}$-monoid objects.
We denote by
${}_*\mathrm{Mon}_{\mathrm{BM}}(\mathcal{C})$
and 
${}_*\mathrm{Mon}_{\mathrm{BM}}^{\rm gen}(\mathcal{C})$
the full subcategories
of $\mathrm{Mon}_{\mathrm{BM}}(\mathcal{C})$
and 
$\mathrm{Mon}_{\mathrm{BM}}^{\rm gen}(\mathcal{C})$
spanned by those functors 
such that the restriction to
$\mathbf{Ass}_L^{\otimes}$ 
has values in final objects of $\mathcal{C}$,
respectively.
The inclusion map
$\mathrm{RM}\to\mathrm{BM}$
induces functors
${}_*\mathrm{Mon}_{\mathrm{BM}}(\mathcal{C})\to
 \mathrm{Mon}_{\mathrm{RM}}(\mathcal{C})$
and
${}_*\mathrm{Mon}_{\mathrm{BM}}^{\rm gen}(\mathcal{C})\to
\mathrm{Mon}_{\mathrm{RM}}^{\rm gen}(\mathcal{C})$
by restriction.
%By restrictionThe inclusion map $\mathrm{RM}\to\mathrm{BM}$
%induces a functor
%${}_*\mathrm{BM}^{\rm gen}(\mathcal{C})\to
%\mathrm{RM}^{\rm gen}(\mathcal{C})$.
We will show that
they give equivalences of $\infty$-categories.

\begin{lemma}\label{lemma:equiv-*-BM-RM}
The functors
${}_*\mathrm{Mon}_{\mathrm{BM}}(\mathcal{C})\to
\mathrm{Mon}_{\mathrm{RM}}(\mathcal{C})$
and
${}_*\mathrm{Mon}_{\mathrm{BM}}^{\rm gen}(\mathcal{C})\to
\mathrm{Mon}_{\mathrm{RM}}^{\rm gen}(\mathcal{C})$
are equivalences
of $\infty$-categories.
\end{lemma}

\proof
We obtain a functor
$\mathrm{Fun}(\mathrm{RM},\mathcal{C})\to
\mathrm{Fun}(\mathrm{BM},\mathcal{C})$
by the right Kan extension
along the inclusion 
$\mathrm{RM}\to \mathrm{BM}$,
which restricts to functors
$\mathrm{Mon}_{\mathrm{RM}}(\mathcal{C})\to
{}_*\mathrm{Mon}{\mathrm{BM}}(\mathcal{C})$
and
$\mathrm{Mon}_{\mathrm{RM}}^{\rm gen}(\mathcal{C})\to
{}_*\mathrm{Mon}_{\mathrm{BM}}^{\rm gen}(\mathcal{C})$.
We can verify that
they are inverses
of ${}_*\mathrm{Mon}_{\mathrm{BM}}(\mathcal{C})
\to \mathrm{Mon}_{\mathrm{RM}}(\mathcal{C})$
and 
${}_*\mathrm{Mon}_{\mathrm{BM}}^{\rm gen}(\mathcal{C})
\to \mathrm{Mon}_{\mathrm{RM}}^{\rm gen}(\mathcal{C})$,
respectively.
\qed

%\newpage

\subsection{Looping construction of right modules}

Let $\mathcal{C}$ be an $\infty$-category
with finite limits.
We consider a generalized
$\mathrm{RM}$-monoid object
$M: \mathrm{RM}\to \mathcal{C}$.
For a morphism $A: *\to M(1)$,
we set
\[ N(0,1)=M(0,1)\times_{M(1)}\{A\},\qquad
   N(1,1)=\{A\}\times_{M(1)}M(1,1)
   \times_{M(1)}\{A\}.\]

The following proposition
can be proved in a similar manner
to Proposition~\ref{prop:construction-bimodule}.

\begin{proposition}\label{prop:construction-right-module}
Let $\mathcal{C}$ be an $\infty$-category with finite limits,
and let $M: \mathrm{RM}\to\mathcal{C}$
be a generalized $\mathrm{RM}$-monoid object of $\mathcal{C}$.
For any morphism $A: *\to M(1)$,
there is an $\mathrm{RM}$-monoid object
$M_A\in \mathrm{Mon}_{\mathrm{RM}}(\mathcal{C})$
such that
$M_A(0,1)\simeq N(0,1)$
and
$M_A(1,1)\simeq N(1,1)$.
\end{proposition}

If $M\in \mathrm{Mon}_{\mathrm{RM}}^{\rm gen}(\mons)$,
then 
we can consider the $\infty$-category
${\rm Alg}_{\mathcal{O}}(M(1))$
of $\mathcal{O}$-algebra objects in $M(1)\in\mons$.
The following theorem can be proved
in a similar manner to
Theorem~\ref{thm:construction-bimodule-O-monoidal}.

\begin{theorem}
\label{thm:construction-right-module-O-monoidal}
Let $M\in \mathrm{Mon}_{\mathrm{RM}}^{\rm gen}(\mons)$
be a generalized $\mathrm{RM}$-monoid
object of $\mons$,
and let
$A\in {\rm Alg}_{\mathcal{O}}(M(1))$.
If 
$M(0,1)^{\otimes}\to M(1)^{\otimes}$ and 
$M(1,1)^{\otimes}\to M(1)^{\otimes}
\times_{\mathcal{O}^{\otimes}}M(1)^{\otimes}$
are cocartesian fibrations of $\infty$-operads,
then $M_A$ is an $\mathrm{RM}$-monoid object
in $\monlax$.
%%such that
%$M_A(0,1)\simeq N(0,1)$
%and $M_A(1,1)\simeq N(1,1)$.
\end{theorem}

%\newpage
%\input{multi-comodule}

\subsection{Multiplicative structures
on right (co)modules}

In this subsection
we study multiplicative structures
on right (co)modules.
We will construct an $\infty$-category
of right modules
%where the ambient $\infty$-category
%is a right module object
in a right module object of $\monoplax$,
and show that it has the structure of an
$\mathcal{O}$-monoidal $\infty$-category.
By taking opposites,
we also show that an $\infty$-category
of right comodules has the structure
of an $\mathcal{O}$-monoidal
$\infty$-category.

We suppose that
$X$
is an $\mathrm{RM}$-monoid
object of $\monoplax$.
We set $X_M=X(0,1)$.
By Remarks~\ref{remark:equivalence-mixed-duoidal-bilax}
and \ref{remark:two-RM-monoid-objects},
we regard $X$ as a mixed
$(\mathbf{RM},\mathcal{O})$-duoidal
$\infty$-category
\[ X
   \longrightarrow \mathbf{RM}^{\otimes}
   \times \mathcal{O}^{\boxtimes,{\rm op}}. \]
As in the case of a mixed $(\mathbf{BM},\mathcal{O})$-duoidal
$\infty$-category in \S\ref{subsection:multiplicative-bimodules},
we obtain a map of opposite $\infty$-operads
\[ \pi: \mathrm{RMod}(X_M)^{\boxtimes,\vee}
    \longrightarrow  \mathrm{Alg}(X_R)^{\boxtimes,\vee} \]
over $\mathcal{O}^{\boxtimes, \rm op}$,
where $\mathrm{RMod}(X_M)^{\boxtimes,\vee}=
\mathrm{Alg}_{\mathrm{RM}/\mathbf{RM}}(X)^{\boxtimes,\vee}$
and
$\mathrm{Alg}(X_R)^{\boxtimes,\vee}=
\mathrm{Alg}_{\mathbf{Ass}/\mathbf{RM}}(X)^{\boxtimes,\vee}$.

We have the following proposition
which is an analogue of Proposition~\ref{prop:bmod-alg-cartesian}.

\begin{proposition}\label{prop:rmod-alg-cartesian}
The map
$\pi: {\rm RMod}(X_M)^{\boxtimes,\vee}\to
{\rm Alg}(X_R)^{\boxtimes,\vee}$
is a cartesian fibration of opposite $\infty$-operads.
\end{proposition}

\proof
We regard $X$ as an $\mathbf{RM}$-monoid
object of $\monoplax$
by Remark~\ref{remark:two-RM-monoid-objects}.
Let $\overline{X}$ be a $\mathbf{BM}$-monoid object
of $\monoplax$ obtained 
by the right Kan extension of $X$
along the inclusion
$\mathbf{RM}^{\otimes}\to\mathbf{BM}^{\otimes}$.
We notice that 
the restriction of $\overline{X}$
to $\mathbf{RM}^{\otimes}$ is equivalent to $X$,
and that 
the restriction 
to $\mathbf{Ass}_L^{\otimes}$
is a final object of $\mathrm{Mon}_{\mathbf{Ass}}(\monoplax)$.
%We regard $Y$ as a mixed
%$(\mathbf{BM},\mathcal{O})$-duoidal
%$\infty$-category
%$Y\to \mathbf{BM}^{\otimes}
%   \times \mathcal{O}^{\boxtimes,{\rm op}}$.
By Proposition~\ref{prop:bmod-alg-cartesian},
we have a cartesian fibration of
opposite $\infty$-operads
$\mathrm{BMod}(\overline{X}_M)^{\boxtimes,\vee}\to
\mathrm{Alg}(\overline{X}_L)^{\boxtimes,\vee}\times_{\mathcal{O}^{\boxtimes,\rm op}}
\mathrm{Alg}(\overline{X}_R)^{\boxtimes,\vee}$.
The proposition follows
by observing that
there are equivalences
$\mathrm{BMod}(\overline{X}_M)^{\boxtimes,\vee}\simeq
\mathrm{RMod}(X_M)^{\boxtimes,\vee}$ and 
$\mathrm{Alg}(\overline{X}_L)^{\boxtimes,\vee}\simeq
\mathcal{O}^{\boxtimes,\rm op}$.
\qed

%@\bigskip

Let $\Gamma\in \mathrm{Bialg}_{(\mathbf{Ass},\mathcal{O})}(X_R)$
be an $(\mathbf{Ass},\mathcal{O})$-bialgebra object of
the mixed $(\mathbf{Ass},\mathcal{O})$-duoidal
$\infty$-category $X_R$.
This determines a map of opposite $\infty$-operads
%an oplax $\mathcal{O}$-monoidal functor
$\Gamma: \mathcal{O}^{\boxtimes,{\rm op}}\to
           {\rm Alg}(X_R)^{\boxtimes,\vee}$
over $\mathcal{O}^{\boxtimes,{\rm op}}$.
We consider the following pullback diagram
\[ \xymatrix{
    {\rm RMod}_{\Gamma}(X_M)^{\boxtimes,\vee}\ar[r]\ar[d]&
    {\rm RMod}(X_M)^{\boxtimes,\vee}\ar[d]^{\pi}\\
    \mathcal{O}^{\boxtimes,{\rm op}}\ar[r]^{\Gamma}&
    {\rm Alg}(X_R)^{\boxtimes,\vee}.\\
}\]
The following theorem can be proved
in a similar manner to Theorem~\ref{thm:Gamma-Bmod-O-monoidal}.

\begin{theorem}\label{thm:Rmod-O-monoidal}
The map
${\rm RMod}_{\Gamma}(X_M)^{\boxtimes,\vee}\to
\mathcal{O}^{\boxtimes,{\rm op}}$ is a cartesian fibration
of opposite $\infty$-operads,
and hence it gives
${\rm RMod}_{\Gamma}(X_M)$
the structure of an $\mathcal{O}$-monoidal $\infty$-category.
Furthermore,
the forgetful functor
${\rm RMod}_{\Gamma}(X_M)\to X_M$
is strong $\mathcal{O}$-monoidal.
\end{theorem}

Next,
we consider a dual situation of the above.
Let $Y$ be an $\mathrm{RM}$-monoid object of $\monlax$.
We regard it as
a mixed $(\mathcal{O},\mathbf{RM})$-duoidal
$\infty$-category
$Y\to \mathcal{O}^{\boxtimes}\times\mathbf{RM}^{\otimes,\rm op}$
by Remark~\ref{remark:equivalence-mixed-duoidal-bilax}.
By Proposition~\ref{prop:duality-mixed-duoidal},
the opposite $Y^{\rm op}\to
\mathbf{RM}^{\otimes}\times\mathcal{O}^{\boxtimes,\rm op}$
has the structure of
a mixed $(\mathbf{RM},\mathcal{O})$-duoidal
$\infty$-category.
By Theorem~\ref{thm:Rmod-O-monoidal},
we obtain a cartesian fibration
$\mathrm{RMod}_{\Gamma}(Y_M^{\rm op})^{\boxtimes,\vee}\to
\mathcal{O}^{\boxtimes,\rm op}$
of opposite $\infty$-operads,
where 
$\Gamma$ is an
$(\mathbf{Ass},\mathcal{O})$-bialgebra object
of the mixed $(\mathbf{Ass},\mathcal{O})$-duoidal
$\infty$-category $Y_R^{\rm op}$.
We can regard $\Gamma$ as an
$(\mathcal{O},\mathbf{Ass})$-bialgebra object
of the mixed $(\mathcal{O},\mathbf{Ass})$-duoidal
$\infty$-category $Y_R$.

\begin{definition}\rm
We define a cocartesian fibration of
$\infty$-operads
\[ \rcomod_{\Gamma}(Y_M)^{\boxtimes}
   \to\mathcal{O}^{\boxtimes} \] 
to be the opposite of
the map
$\mathrm{RMod}_{\Gamma}(Y_M^{\rm op})^{\boxtimes,\vee}\to
\mathcal{O}^{\boxtimes}$.
We call 
$\rcomod_{\Gamma}(Y_M)$
the $\infty$-category of 
right $\Gamma$-comodules of $Y_M$.
\end{definition}

Roughly speaking,
an object of
$\rcomod_{\Gamma}(Y_M)$
is an object $M_x$ of $(Y_M)_x$ for some $x\in\mathcal{O}$
equipped with a coaction map
\[ M_x\longrightarrow M_x\otimes_{x} \Gamma_{x} \]
that is coassociative and counital
up to higher coherent homotopy.   

In summary we obtain the following corollary.

\begin{corollary}\label{cor:mixed-to-product-right-comodules}
Let $Y$ be an $\mathrm{RM}$-monoid object of $\monlax$. 
For an $(\mathcal{O},\mathbf{Ass})$-bialgebra object
$\Gamma$ 
of the mixed $(\mathcal{O},\mathbf{Ass})$-duoidal
$\infty$-category $Y_R$,
the $\infty$-category $\rcomod_{\Gamma}(Y_M)$
of right $\Gamma$-comodules in $Y_M$
has the structure of an $\mathcal{O}$-monoidal
$\infty$-category such that
the forgetful functor
$\rcomod_{\Gamma}(Y_M)\to Y_M$
is strong $\mathcal{O}$-monoidal.
\end{corollary}

%\newpage

\subsection{Right comodules in
$(\mathcal{O}\times\mathbf{Ass})$-monoidal
$\infty$-categories}

In this subsection
we start
with an $(\mathcal{O}\times\mathbf{Ass})$-monoidal
$\infty$-category $\mathcal{M}$ and
an $(\mathcal{O}\times\mathbf{Ass})$-algebra $A$ in $\mathcal{M}$,
and
construct a $\mathrm{RM}$-monoid object
$\mathbb{RM}_A(\mathcal{M})$ in $\monlax$
such that $\mathbb{RM}_A(\mathcal{M})(0,1)
\simeq \mathrm{RMod}_A(\mathcal{M})$.
%By applying Theorem~\ref{thm:construction-right-module-O-monoidal}
By applying Corollary~\ref{cor:mixed-to-product-right-comodules}
to $\mathbb{RM}_A(\mathcal{M})$, 
we obtain an $\mathcal{O}$-monoidal structure
on the $\infty$-category of
right comodules in $\mathcal{M}$.

We apply Theorem~\ref{thm:construction-right-module-O-monoidal}
to the following situation:
We recall that
for an $(\mathcal{O}\times\mathbf{Ass})$-monoid object
$\mathcal{M}$ of $\cat^{\rm geo}$,
we can construct a generalized $\mathrm{BM}$-monoid
object $\mathbb{BM}(\mathcal{M})$
of $\mathrm{Mon}_{\mathcal{O}}^{\rm geo}$
such that $\mathbb{BM}(\mathcal{M})(i,j)
\simeq \mathrm{BMod}(\mathcal{M})$
and $\mathbb{BM}(\mathcal{M})(i)\simeq
\mathrm{Alg}(\mathcal{M})$.
%by \cite[Lemma~4.19 and Theorem~4.39]{Haugseng1}.
Let $\mathrm{RM}'$ be a full subcategory
of $\mathrm{RM}$ spanned by
all objects except $(0)$.
We denote by $i: \mathrm{RM}'\to\mathrm{BM}$
the inclusion functor.
We set ${}_*\mathbb{BM}(\mathcal{M})=
i_*i^*\mathbb{BM}(\mathcal{M})$,
where $i^*: \mathrm{Fun}(\mathrm{BM},\mons^{\rm geo})
\to \mathrm{Fun}(\mathrm{RM}',\mons^{\rm geo})$
is the restriction,
and $i_*: \mathrm{Fun}(\mathrm{RM}',\mons^{\rm geo})
\to \mathrm{Fun}(\mathrm{BM},\mons^{\rm geo})$
is its right adjoint given
by the right Kan extension along $i$.
We can verify that ${}_*\mathbb{BM}(\mathcal{M})$
is an object of
${}_*\mathrm{Mon}_{\mathrm{BM}}^{\rm gen}(\mons^{\rm geo})$.
By Lemma~\ref{lemma:equiv-*-BM-RM},
we obtain a generalized $\mathrm{RM}$-monoid
object $\mathbb{RM}(\mathcal{M})$ of $\mons^{\rm geo}$
such that $\mathbb{RM}(\mathcal{M})(0,1)\simeq
{\rm RMod}(\mathcal{M})$,
$\mathbb{RM}(\mathcal{M})(1,1)\simeq {\rm BMod}(\mathcal{M})$,
and $\mathbb{RM}(\mathcal{M})(1)\simeq {\rm Alg}(\mathcal{M})$.

\if0
We denote by 
$p: {\rm RMod}(\mathcal{M})^{\boxtimes}\to
   \mathcal{O}^{\boxtimes}$
a cocartesian fibration of $\infty$-operads
which gives ${\rm RMod}(\mathcal{M})$
the structure of the $\mathcal{O}$-monoidal $\infty$-category.
%We also denote by 
%$q: {\rm Alg}(\mathcal{M})^{\boxtimes}\to
%   \mathcal{O}^{\boxtimes}$
%a cocartesian fibration of $\infty$-operads
%which gives ${\rm Alg}(\mathcal{M})$
%the structure of the $\mathcal{O}$-monoidal
%$\infty$-category.
We have a commutative diagram
\[ \xymatrix{
  {\rm RMod}(\mathcal{M})^{\boxtimes}
  \ar[rr]^-{\pi}\ar[dr]_p&&
     {\rm Alg}(\mathcal{M})^{\boxtimes}
           \ar[dl]^q\\
   &\mathcal{O}^{\boxtimes}.&\\
  }\]
\fi

The following proposition
can be proved
in a similar manner to
Proposition~\ref{prop:construction-BMod-object}.
%and
%Lemma~\ref{lemma:BMod-Alg-cocartesian},
%respectively.

\begin{proposition}\label{prop:construction-RMod-object}
Let $\mathcal{M}\in {\rm Mon}_{(\mathcal{O}\times\mathbf{Ass})}
(\cat^{\rm geo})$.
For any $(\mathcal{O}\times\mathbf{Ass})$-algebra
$A$ of $\mathcal{M}$,
there is an $\mathrm{RM}$-monoid object
$\mathbb{RM}_A(\mathcal{M})$ in $\mons^{\rm lax, geo}$
such that 
$\mathbb{RM}_A(\mathcal{M})(0,1)\simeq
{\rm RMod}_A(\mathcal{M})$
and
$\mathbb{RM}_A(\mathcal{M})(1,1)
\simeq {}_A{\rm BMod}_A(\mathcal{M})$.
%in $\mons^{\rm lax, geo}$ via relative tensor products.   
\end{proposition}

We can regard the $\mathrm{RM}$-monoid
object $\mathbb{RM}_A(\mathcal{M})$ in $\mons^{\rm lax, geo}$
as a mixed $(\mathcal{O},\mathbf{RM})$-duoidal
$\infty$-category
by Remark~\ref{remark:equivalence-mixed-duoidal-bilax}.
Let $\Gamma$ 
be an $(\mathcal{O},\mathbf{Ass})$-bialgebra objects
of the $(\mathcal{O},\mathbf{Ass})$-duoidal
$\infty$-categories $\mathbb{RM}_A(\mathcal{M})_R$
We notice that $\Gamma$ 
corresponds to an $(\mathcal{O},\mathbf{Ass})$-bialgebra
of ${}_A\mathrm{BMod}_A(\mathcal{M})$.

\begin{definition}\rm
We write 
\[ \rcomod_{(A,\Gamma)}
   (\mathcal{M})^{\boxtimes}
   \longrightarrow \mathcal{O}^{\boxtimes} \]
for the cocartesian fibration of
$\infty$-operads
$\mathrm{RMod}_{\Gamma}
(\mathbb{RM}_A(\mathcal{M})_M)^{\boxtimes}
\to \mathcal{O}^{\boxtimes}$,
and call it
%the underlying $\infty$-category
%$\rcomod_{(A,\Gamma)}(\mathcal{M})$
the $\infty$-category of 
right $(A,\Gamma)$-comodules of $\mathcal{M}$.
\end{definition}

\if0
By applying Corollary~\ref{cor:mixed-to-product-comodules}
to ${}_A\mathbb{BM}_B(\mathcal{M})$,
we obtain the following theorem.

Let $\Gamma$ be an $(\mathbf{Ass},\mathcal{O})$-bialgebra
of the mixed $(\mathbf{Ass},\mathcal{O})$-duoidal
$\infty$-category $\mathbb{RM}_A(\mathcal{M})_R$.
We notice that
it corresponds to an $(\mathcal{O},\mathbf{Ass})$-bialgebra
object of ${}_A{\rm BMod}_A(\mathcal{M})$.
By applying Theorem~\ref{thm:Rmod-O-monoidal},
we obtain a
${\rm RMod}_{\Gamma}(\mathbb{RM}_A(\mathcal{M})_M)$
such that the forgetful functor
${\rm RMod}_{\Gamma}({\rm RMod}_A(\mathcal{M})^{\rm op})
\to {\rm RMod}_A(\mathcal{M})^{\rm op}$
is strong $\mathcal{O}$-monoidal.

\begin{definition}\rm
We write
\[ \rcomod_{(A,\Gamma)}(\mathcal{M})^{\boxtimes}
   \longrightarrow \mathcal{O}^{\boxtimes} \]
for the map
$\rcomod_{\Gamma}(\mathbb{RM}_A(\mathcal{M})_M)^{\boxtimes}
   \to \mathcal{O}^{\boxtimes}$.
We call $\rcomod_{(A,\Gamma)}(\mathcal{M})$
the $\infty$-category of right
$(A,\Gamma)$-comodules of $\mathcal{M}$.
\end{definition}
\fi

In summary,
we obtain the following theorem.

\begin{theorem}
%[{cf.~Theorem~\ref{thm:main-theorem}}]
\label{thm:right-comodule-main-theorem}
Let $\mathcal{M}$ be an $(\mathcal{O}\times\mathbf{Ass})$-monoid
object of $\cat^{\rm geo}$, and
let $A$ be an $(\mathcal{O}\times\mathbf{Ass})$-algebra
object of $\mathcal{M}$.
For an $(\mathcal{O}, \mathbf{Ass})$-bialgebra
$\Gamma$ in ${}_A{\rm BMod}_A(\mathcal{M})$,
the $\infty$-category
$\rcomod_{(A,\Gamma)}(\mathcal{M})$
of right $(A,\Gamma)$-comodules of $\mathcal{M}$
has the structure of an $\mathcal{O}$-monoidal
$\infty$-category such that 
the forgetful functor
$\rcomod_{(A,\Gamma)}(\mathcal{M})\to
{\rm RMod}_A(\mathcal{M})$
is strong $\mathcal{O}$-monoidal.
\end{theorem}

\begin{corollary}
Let $\mathcal{M}$ be an $\mathbf{Ass}^{n+1}$-monoid
object in $\cat^{\rm geo}$
and
let $A$ be an $\mathbf{Ass}^{n+1}$-algebra
object in $\mathcal{M}$.
For an $(\mathbf{Ass}^n, \mathbf{Ass})$-bialgebra
$\Gamma$ in ${}_A{\rm BMod}_A(\mathcal{M})$,
the $\infty$-category
$\rcomod_{(A,\Gamma)}(\mathcal{M})$
of right $(A,\Gamma)$-comodules of $\mathcal{M}$
has the structure of an $\mathbf{Ass}^n$-monoidal
$\infty$-category such that
the forgetful functor
$\rcomod_{(A,\Gamma)}(\mathcal{M})\to
{\rm RMod}_A(\mathcal{M})$
is strong $\mathbf{Ass}^n$-monoidal.
\end{corollary}

%\newpage
%{\footnotesize
%\input{ref}

%}


\begin{thebibliography}{99}

\bibitem{Adams1}
J. F. Adams, 
On the structure and applications of the Steenrod algebra,
Comment. Math. Helv. 32 (1958), 180--214.
  
\bibitem{Adams2}  
J. F. Adams, 
On the non-existence of elements of Hopf invariant one,
Ann. of Math. (2) 72 (1960), 20--104.

\bibitem{Aguiar-Mahajan}
M. Aguiar and S. Mahajan, 
Monoidal functors, species and Hopf algebras,
CRM Monograph Series, 29. 
American Mathematical Society, Providence, RI, 2010. 

\bibitem{Barwick}
C. Barwick, 
From operator categories to higher operads,
Geom. Topol. 22 (2018), no. 4, 1893--1959. 

\bibitem{CH}
H. Chu and R. Haugseng, 
Homotopy-coherent algebra via Segal conditions,
Adv. Math. 385 (2021), Paper No. 107733, 95 pp.

\bibitem{Gepner-Haugseng}
D. Gepner and R. Haugseng, 
Enriched $\infty$-categories via non-symmetric $\infty$-operads. 
Adv. Math. 279 (2015), 575--716. 

%\bibitem{GHN}
%D. Gepner, R. Haugseng, and T. Nikolaus, 
%Lax colimits and free fibrations in ∞-categories.
%Doc. Math. 22 (2017), 1225--1266. 

%\bibitem{Glasman}
%S. Glasman, 
%A spectrum-level Hodge filtration on 
%topological Hochschild homology. 
%Selecta Math. (N.S.) 22 (2016), no. 3, 1583--1612. 

\bibitem{Haugseng1}
R. Haugseng, 
The higher Morita category of $\mathbb{E}_n$-algebras,
Geom. Topol. 21 (2017), no. 3, 1631--1730. 

\bibitem{HHLN1}
R. Haugseng, F. Hebestreit, S. Linskens, and J. Nuiten,
Lax monoidal adjunctions, two-variable fibrations and 
the calculus of mates,
Proc. Lond. Math. Soc. (3) 127 (2023), no.4, 889--957.

\bibitem{Lurie1}
J. Lurie, 
Higher topos theory,
Annals of Mathematics Studies, 170. 
Princeton University Press, Princeton, NJ, 2009.

\bibitem{Lurie2}
J. Lurie,
Higher algebra,
%(version 9/14/2014),
available at 
http://www.math.harvard.edu/\~{}lurie/.

%\bibitem{Milnor}
%J. Milnor, 
%The Steenrod algebra and its dual,
%Ann. of Math. (2) 67 (1958), 150--171.

\bibitem{MRW}
H. R. Miller, D. C. Ravenel, and W. S. Wilson, 
Periodic phenomena in the Adams-Novikov spectral sequence,
Ann. of Math. (2) 106 (1977), no. 3, 469--516. 

\bibitem{Morava}
J. Morava, 
Noetherian localisations of categories of cobordism comodules,
Ann. of Math. (2) 121 (1985), no. 1, 1--39. 

\bibitem{Ravenel}
D. C. Ravenel, 
Complex cobordism and stable homotopy groups of spheres,
Pure and Applied Mathematics, 121. 
Academic Press, Inc., Orlando, FL, 1986. 

\bibitem{Stefanich}
G. Stefanich,
Higher sheaf theory I: Correspondences,
preprint,
arXiv:2011.03027.
  
\bibitem{Street}
R. Street, 
Monoidal categories in, and linking, geometry and algebra,
Bull. Belg. Math. Soc. Simon Stevin 19 (2012), no. 5, 769--821. 

\bibitem{Torii}
T. Torii,
On Quasi-Categories of Comodules and Landweber Exactness,
In: Ohsawa T., Minami N. (eds) 
Bousfield Classes and Ohkawa's Theorem. BouCla 2015. 
Springer Proceedings in Mathematics \& Statistics, vol~309,  
325--380. Springer, Singapore, 2020.

%\bibitem{Torii6}
%T. Torii,
%A perfect pairing for monoidal adjunctions,
%Proc. Amer. Math. Soc. 151 (2023), no.12, 5069--5080.

%\bibitem{Torii7}
%T. Torii,
%Uniqueness of monoidal adjunctions,
%Homology, Homotopy and Applications,
%vol. 26 (2), 2024, 259--272
  
\bibitem{Torii1}
T. Torii,
On duoidal $\infty$-categories, 
%preprint,
%arXiv:2106.14121,
J. Homotopy Relat. Struct. 20 (2025), no. 1, 125--162.
    
\bibitem{Torii2}
T. Torii,
On higher monoidal $\infty$-categories,
preprint,
arXiv:2111.00158.

\bibitem{Torii4}
T. Torii,
Duoidal $\infty$-categories of operadic modules,
preprint,
arXiv:2204.11152.
  
\bibitem{Torii5}
T. Torii,
Map monoidales and duoidal $\infty$-categories,
preprint,
arXiv:2406.00223.
  
\end{thebibliography}
\end{document}